%% file: article.tex
\documentclass[preprint,sort&compress]{elsarticle}
\usepackage[margin=1in]{geometry}
\usepackage{lineno}
\usepackage{amssymb,amsxtra,amsmath,amsthm}
\usepackage{mathrsfs}
\usepackage{latexsym}
\usepackage{ifthen}
\usepackage{graphicx}
\usepackage{calc,fancyhdr}
\usepackage{indentfirst}
\usepackage{parskip,algorithm,algorithmicx}
\usepackage{algpseudocode}
\usepackage{color}
\usepackage{bigints}
\usepackage{subfig}
\usepackage{booktabs}

\usepackage[bottom]{footmisc}

\usepackage{Article}
\usepackage{color}
\newcommand{\revone}[1]{\color{black}#1\normalcolor}
\theoremstyle{remark}
\newtheorem{remark}{Remark}

\journal{JCP}

\begin{document}

\begin{frontmatter}

\title{An Efficient High-Order Meshless Method for Advection-Diffusion Equations on Time-Varying Irregular Domains}

\author[addr1]{Varun Shankar\corref{corresp}}
\address[addr1]{School of Computing, University of Utah, UT, USA}
\ead{shankar@cs.utah.edu}
\cortext[corresp]{Corresponding Author}

\author[addr2]{Grady B. Wright}
\address[addr2]{Department of Mathematics, Boise State University, ID, USA}
\ead{gradywright@boisestate.edu}

\author[addr3]{Aaron L. Fogelson}
\address[addr3]{Departments of Mathematics and Biomedical Engineering, University of Utah, UT, USA}
\ead{fogelson@math.utah.edu}

\begin{abstract}
We present a high-order radial basis function finite difference (RBF-FD) framework for the solution of advection-diffusion equations on time-varying domains. Our framework is based on a generalization of the recently developed Overlapped RBF-FD method that utilizes a novel automatic procedure for computing RBF-FD weights on stencils in variable-sized regions around stencil centers. This procedure eliminates the overlap parameter $\delta$, thereby enabling tuning-free assembly of RBF-FD differentiation matrices on moving domains. In addition, our framework utilizes a simple and efficient procedure for \emph{updating} differentiation matrices on moving domains tiled by node sets of time-varying cardinality. Finally, advection-diffusion in time-varying domains is handled through a combination of rapid node set modification, a new high-order semi-Lagrangian method that utilizes the new tuning-free overlapped RBF-FD method, and a high-order time-integration method. The resulting framework has no tuning parameters and has $O(N \log N)$ time complexity. We demonstrate high-orders of convergence for advection-diffusion equations on time-varying 2D and 3D domains for both small and large Peclet numbers. We also present timings that verify our complexity estimates. Finally, we utilize our method to solve a coupled 3D problem motivated by models of platelet aggregation and coagulation, once again demonstrating high-order convergence rates on a moving domain.
\end{abstract}

\end{frontmatter}

\input{Intro}
\input{OverlapNew}
\input{Methods}
\input{LinAlg}
\input{CompAnalysis}

\input{Results}

\input{coupled}
\input{Discussion}

\section*{Acknowledgments}
VS was supported by NSF grants CISE CCF 1714844 and DMS-1521748. ALF was supported by NSF grant DMS-1521748 and NIHBL grant 1U01HL143336. GBW was supported by NSF grants CISE CCF 1717556 and DMS 1952674.


\section*{References}
\bibliography{article_refs_mod}

\end{document}

%% file: Intro.tex
\section{Introduction}
\label{sec:intro}

Collocation methods based on radial basis functions (RBFs) have been increasingly popular for numerically solving partial differential equations (PDEs), due to their high-order convergence rates and their ability to naturally handle scattered node layouts on arbitrary domains. RBF interpolants can be used to generate both pseudospectral (RBF-PS) and finite-difference (RBF-FD) methods~\cite{Bayona2010,Davydov2011,Wright200699,FlyerNS,FlyerPHS,BarnettPHS}. RBF-based methods are also easily applied to the solution of PDEs on node sets that are not unisolvent for polynomials, such as ones lying on the sphere $\mathbb{S}^2$~\cite{FlyerWright:2007,FlyerWright:2009,FoL11,FlyerLehto2012} and other general surfaces~\cite{Piret2012, Piret2016,FuselierWright2013,SWFKJSC2014,LSWSISC2017}.

\revone{The} focus of this paper is on advection-diffusion problems on domains $\Omega(t)$ with boundary conditions enforced at time-varying internal embedded boundaries and a fixed outer boundary. This can be modeled by the following equations:
\begin{align}
\frac{\partial c}{\partial t} + \vu \cdot \nabla c &= \nu \Delta c + f(\vx,t), \vx \in \Omega(t),\label{eq:adv-diff} \\
\alpha(\vx,t) \vn \cdot \nabla c + \beta(\vx,t) c &= g(\vx,t), \vx \in \partial \Omega(t) \label{eq:ad-bc},
\end{align}
where $c(\vx,t)$ is a scalar quantity being transported in the incompressible velocity field $\vu(\vx,t)$, $\nu$ is the diffusion coefficient, $\alpha$ and $\beta$ are functions that determine boundary conditions (linear in this article), $\vn$ are the unit \emph{outward} normals to the domain, and $g(\vx,t)$ is either a prescribed or numerically computed boundary condition. Our interest in the above equations stems from their application in the modeling and simulation of platelet aggregation and coagulation~\cite{LEIDERMAN:2011:GWF,LEIDERMAN:2013:IHT,LEIDERMAN:2014:OMM}. Broadly speaking, numerical methods to solve such systems can be divided into three categories: (a) \textbf{Eulerian} methods (b) \textbf{Lagrangian} methods and (c) \textbf{semi-Lagrangian} (SL) methods.

There is extensive literature on Eulerian finite difference (FD) or finite volume (FV) methods for solving PDEs along with boundary conditions on fixed irregular surfaces embedded in the computational domain.  Such methods (most of which are designed for Cartesian grids) are mainly of two types. The first type involves \emph{augmenting} the FD/FV scheme to enforce boundary conditions at the irregular boundary. This could be done via \emph{spreading} and \emph{restriction} as in the immersed boundary (IB) method~\cite{IBM1,IBM2,PESKIN:2002:IBM}, or via adding unknowns to the system to force the PDE to satisfy boundary conditions (at the irregular embedded boundaries) as in the wide class of \emph{forcing methods}~\cite{IB-externalforce-Goldstein-JCP1993,ForcingBC-Mohd-Yusof1997,IB-forcing_sink-KimChoi-JCP2001,YaoFogelson2012,SWFKIJNMF2014}, the ghost cell method~\cite{ColoniusTaira08,Glowinski}, and the more recent immersed boundary smooth extension method~\cite{IBSE1,IBSE2}. In contrast, the second type involves modifying FD/FV stencils near the boundary, such as in the original direct forcing method~\cite{IB-forcing-Fadlun-JCP2000}, the immersed interface method (IIM)~\cite{IIM1}, the embedded boundary method (EBM)~\cite{EB-poisson-JCP1998}, the sharp interface method~\cite{Udaykumar1999,IB-cutcell-Ye-JCP1999}, and the capacity function finite volume method~\cite{Calhoun99acartesian}. Finally, to tackle moving boundaries outside the original IB framework, a common approach involves \emph{converting} the moving-boundary problem into a series of fixed-boundary problems each solved by one of the above approaches (\emph{e.g.}, see~\cite{McCorquodaleColellaJohansen01}). Such methods typically require an additional spatial extrapolation step to fill newly-uncovered grid points as the domain boundary moves. For all these types of Eulerian methods, obtaining a stable, high-order discretization in space and time can be challenging both due to the presence of a background Cartesian grid and the need for spatial extrapolation to fill newly-uncovered grid points.

In contrast to Eulerian methods that use a fixed background grid, Lagrangian methods involve populating the moving domain with a set of marker particles that move with the velocity field $\vu$. In this case, the advection term is handled without any difficulty. However, to discretize the diffusion term, one of the following approaches can be used: (a) interpolate quantities to a fixed background grid and discretize the PDE there like in the material point method (MPM)~\cite{Gritton2017}; (b) discretize the diffusion term directly on the distorted Lagrangian grid as in smoothed particle hydrodynamics (SPH)~\cite{trask2015Scalable}); or (c) \revone{some variation on the weighted particle method~\cite{PSE1} (sometimes called the particle strength exchange method)}. In all cases, it may be necessary to introduce some form of Lagrangian particle rearrangement to improve spatial resolution and convergence rates~\cite{LagRe}. Alternatively, it is possible to \revone{reformulate the PDE using the Feynman-Kac formula} so that the diffusion term is also handled in a Lagrangian fashion~\cite{FFSL,BVSDE1,BVSDE2}.

Semi-Lagrangian (SL) methods are Eulerian methods that use fictitious Lagrangian particles to determine the numerical domain of dependence. We focus on the class of \emph{backward} SL methods, which have found wide application to problems in fluid dynamics, climate modeling, and numerical weather prediction~\cite{staniforth1991SL,StaniforthWoodJCP2008,SmolPudy92,SmolarkiewiczMargolin97,Xiu2002}. For a pure advection equation, these methods assume that Lagrangian marker particles have arrived at every time-step on an Eulerian grid (or more generally node set). By tracing these particles backward through the velocity field (and to a previous time level), determining their departure positions, and interpolating the solution to those departure positions from the fixed Eulerian node set, one can determine how much material was advected to a given Eulerian location. Solving an advection-diffusion problem then amounts to using an appropriate \revone{splitting} scheme. The advantage of this method is that the diffusion operator is always discretized on the Eulerian grid. In addition, when solving problems with moving boundaries, the SL framework obviates the need for any spatial extrapolation.

The numerical method presented in this paper relies on the SL framework for precisely these reasons. While RBF methods have been used within the SL framework before, these either were global RBF methods~\cite{TELA:TELA0009}, methods that relied on Voronoi cells~\cite{iske2002}, or localized RBF methods designed specifically for the sphere~\cite{SWJCP2018}. Our new method is based on a generalization of the overlapped RBF-FD method~\cite{ShankarJCP2017,SFJCP2018,SNKJCP2018,SNWSISC2020}, and therefore allows the use of scattered or quasi-uniform nodes in place of a background Cartesian grid, allowing for arbitrary \revone{outer} boundaries. The time-varying nature of the domain is handled by enabling or disabling pre-existing background nodes contained by the moving boundaries (and in a \emph{small} neighborhood around them). To facilitate this node set adaptation, we represent the moving boundaries using a high-order accurate parametric model built from Lagrangian markers. In this way, our node sets always conform to the time-varying domain, unlike the Eulerian methods discussed above. This technique allows efficient updates to differentiation matrices and necessitates recalculation of overlapped RBF-FD weights \emph{only} in neighborhoods around the moving boundaries. The resulting meshless method allows for high orders of spatial and temporal convergence, does not require spatial extrapolation, and is of \revone{quasi-linear} computational complexity.

The remainder of the paper is organized as follows. In Section 2, \revone{we present our generalization of the overlapped RBF-FD method that removes tuning parameters}. In Section \revone{3}, we present and describe our overall numerical method in Algorithm \ref{alg:sl_master}, complete with error estimates and parameter choices. We present a simple and efficient preconditioner in Section \revone{4}, which we then use to solve the time-varying sparse linear system resulting from our numerical method. We conduct a thorough complexity analysis of Algorithm \ref{alg:sl_master} in Section \revone{5}. Then, in Section \revone{6}, we present 2D and 3D convergence tests on problems with moving embedded boundaries for a range of Peclet numbers. Finally, in Section \revone{7}, we present an application of our method to a 3D coupled problem with time-varying boundary conditions on a moving domain inspired by mathematical models of platelet aggregation and coagulation. We conclude with a summary and comments on future work in Section \revone{8}.

%% file: OverlapNew.tex
\revone{\section{An automatic overlapped RBF-FD method}
\label{sec:new_overlap}
We now present a generalization of the overlapped RBF-FD method that eliminates the overlap parameter $\delta$, and instead automatically computes, tests, and retains/discards candidate weights on a given stencil. Our approach for automation is to use two stability indicators to indicate whether a set of computed weights is of sufficient quality. In the discussion that follows, we will primarily focus on the new method, remarking on the older version presented in~\cite{ShankarJCP2017,SFJCP2018,SNKJCP2018} as needed.  

Let $X = \{\vx_k\}_{k=1}^N$ be a global set of nodes on a domain $\Omega \subset \mathbb{R}^d$. Define the stencil $P_k$ to be the set of nodes containing node $\vx_{\calI^k_1}$ and its $n-1$ nearest neighbors $\{\vx_{\calI^k_2},\hdots,\vx_{\calI^k_n}\}$; here, $\{\calI^k_1,\hdots,\calI^k_n\}$ are indices that map into the global node set $X$ and $\calI^k_1 = k$. Without loss of generality, we focus on the stencil $P_1$. Let $1 \leq p_1 \leq n$ be the number of points on this stencil for which we wish to compute RBF-FD weights. Further, define $R_1$ to be the global indices of these $p_1$ nodes so that
\begin{align}
R_1 = \{\calR^1_1, \calR^1_2, \hdots, \calR^1_{p_1}\}.
\label{eq:ball_inds}
\end{align}
Next, let $\mathbb{B}_1$ be the ball containing the nodes whose indices are in $R_1$. Thus, 
\begin{align}
\mathbb{B}_1 = \{\vx_{\calR^1_1}, \hdots, \vx_{\calR^1_{p_1}} \}.
\label{eq:ball_pts}
\end{align}
We discuss how to obtain $\mathbb{B}_1$ and $\mathbb{R}_1$ in Section \ref{sec:influence}. First, in Section \ref{sec:weights_comp}, we describe how to compute RBF-FD weights for all the nodes in the ball $\mathbb{B}_1$. 

\subsection{Computing weights}
\label{sec:weights_comp}

The weights for all the nodes in $\mathbb{B}_1$ with indices in $R_1$ are computed using the following augmented local RBF interpolant on $P_1$:
\begin{align}
s_1(\vx,\vy) = \sum\limits_{j=1}^n w^1_j(\vy) \|\vx - \vx_{\calI^1_j}\|^m + \sum\limits_{i=1}^{M} \lambda^1_i(\vy) \psi^1_i(\vx),
\label{eq:rbf_interp}
\end{align}
where $\|\vx - \vx_{\calI^1_j}\|^m$ is the polyharmonic spline (PHS) RBF of degree $m$ ($m$ is odd), and $\{\psi^1_i(\vx)\}$ form a basis for the space of polynomials of total degree degree $\ell$ in $d$ dimensions so that $M = {\ell + d \choose d}$; common choices for these include monomials~\cite{SFJCP2018} or orthogonal polynomials~\cite{SNKJCP2018}. In this work, we select the $\psi^1_i(\vx)$ functions to be Legendre polynomials. The $n$ overlapped RBF-FD weights associated with the point $\vy$ are $w^1_j(\vy), j=1,\hdots,n$. We compute the weights for the linear operator $\calL$ uniquely at all nodes in $\mathbb{B}_1$ with indices in the set $R_1$ by imposing the following two (sets of) conditions on \eqref{eq:rbf_interp}:
\begin{align}
\lf.s_1 \rt|_{\vx \in P_1, \vy \in \mathbb{B}_1} &= \lf.\calL \|\vx - \vx_{\calI^1_j}\|^m\rt|_{\vx \in \mathbb{B}_1}, j=1,\hdots,n, \label{eq:interp_constraint}\\
\sum_{j=1}^n \lf.w_j^1(\vy) \psi_i^1(\vx) \rt|_{\vx \in P_1, \vy \in \mathbb{B}_1} &= \lf.\calL \psi^1_i(\vx)\rt|_{\vx \in \mathbb{B}_1}, i=1,\hdots,M. \label{eq:poly_constraint}
\end{align}
These conditions enforce that the weights are exact for both $\calL$ applied to the PHS RBF and to the polynomial basis. In this work, we use the heuristic $n = 2M+1$~\cite{ShankarJCP2017,SFJCP2018,FlyerNS,FlyerPHS}; however, larger stencil sizes can sometimes be beneficial~\cite{BayonaBoundary}. The constraints \eqref{eq:interp_constraint}--\eqref{eq:poly_constraint} for determining the weights in \eqref{eq:rbf_interp} can be collected into the following block linear system:
\begin{align}
\begin{bmatrix}
A_1 & \Psi_1 \\
\Psi_1^T & 0
\end{bmatrix}
\begin{bmatrix}
W_1 \\
W^{\psi}_1
\end{bmatrix}
=
\begin{bmatrix}
B_{A_1} \\
B_{\Psi_1}
\end{bmatrix},
\label{eq:rbf_linsys}
\end{align}
where
\begin{align}
(A_1)_{ij} &= \|\vx_{\calI^1_i} - \vx_{\calI^1_j} \|^m, i,j=1,\hdots,n, \\
(\Psi_1)_{ij} &= \psi^1_j(\vx_{\calI^1_i}), i=1,\hdots,n, j=1,\hdots,M,\\
(B_{A_1})_{ij} &= \lf.\calL \|\vx - \vx_{\calI^1_i} \|^m \rt|_{\vx = \vx_{\calR^1_j}}, i=1,\hdots,n, j=1,\hdots,p_1, \\
(B_{\Psi_1})_{ij} &= \lf.\calL \psi^1_i(\vx)\rt|_{\vx = \vx_{\calR^1_j}}, i=1,\hdots,M, j=1,\hdots,p_1.
\end{align}
$W_1$ is the matrix of overlapped RBF-FD weights, with each column containing the RBF-FD weights for a point $\vx \in \mathbb{B}_1$. It is also notationally useful to refer to the column of $W_1$ associated with the point $\vy$ as $W_1(\vy)$. The linear system \eqref{eq:rbf_linsys} has a unique solution if the nodes in $P_1$ are distinct and $\Psi_1$ has full column rank~\cite{Fasshauer:2007,Wendland:2004}. The matrix of polynomial coefficients $W^{\psi}_1$ enforces the polynomial reproduction constraint \eqref{eq:poly_constraint}. This constraint ensures that the \emph{local} approximation error is bounded by $O \lf(h^{\ell+1-\theta} \rt)$, where $\theta$ is the order of the differential operator $\calL$, and $h$ is the largest distance between the point at which the weights are computed and every other point in the stencil~\cite{DavydovSchaback2018}. 

\subsection{Automatically determining $R_1$ and $\mathbb{B}_1$}
\label{sec:influence}

In previous versions of the overlapped RBF-FD method, the set $R_1$ and the ball $\mathbb{B}_1$ were determined by defining a overlap parameter $\delta \in (0,1]$ such that all nodes within a distance $(1-\delta) \rho_1$ from the center lay within the ball $\mathbb{B}_1$, where $\rho_1 = \max\limits_{1 \leq j \leq n} \|\vx_{\calI^1_1} - \vx_{\calI^1_j}\|$. As such an approach requires tuning $\delta$, in this work, we present an automatic approach to determine the sets $R_1$ and $\mathbb{B}_1$. This new approach uses a pair of stability indicators to determine whether a set of computed weights is of sufficient quality. In the discussion that follows, we continue to use $\vy$ to refer to the point at which weights ${\bf w}^1(\vy)$ are computed, and $\vx$ to refer to points comprising the stencils.

\subsubsection{The local $\calL$-Lebesgue function indicator}
\label{sec:lebesgue}

As described in~\cite{ShankarJCP2017} and noted in~\cite{FlyerElliptic,Bayona2019}, local $\calL$-Lebesgue functions play a key role in assessing the suitability of a set of RBF-FD weights. Large values of the local $\calL$-Lebesgue function can lead to spurious eigenvalues in the differentiation matrix corresponding to $\calL$. This fact was used to develop a stability indicator for the overlapped RBF-FD method in~\cite{ShankarJCP2017} which was used to discard unsuitable weights on a given stencil. Recall that $\vx_{\calI^1_1}$ and its $n-1$ neighbors form the stencil $P_1$. Then, letting $\vy_1 = \vx_{\calI^1_1}$, one can define the $\calL$-Lebesgue function at $\vx_1$ as the $\ell_1$-norm of the weight vector at that point:
\begin{align}
\Lambda_{\calL}(\vy) = \|W_1(\vy)\|_1,
\end{align}
where $W_1(\vy) ={\bf w}^1(\vy)$ refers to the column of $W_1$ corresponding to the point $\vy$. We can now define a set $\mathbb{B}_{11}$ as:
\begin{align}
\mathbb{B}_{11} = \{\vy \in P_1 \mid \Lambda_{\calL}(\vy) \leq \Lambda_{\calL}(\vy_1) \},
\label{eq:leb_set}
\end{align}
\emph{i.e.}, the set of all points in $P_1$ where the $\calL$-Lebesgue function takes on values smaller than at the stencil center $\vy_1$. In~\cite{ShankarJCP2017}, the RBF-FD weights from lower-order methods were computed by solving \eqref{eq:rbf_linsys} (with the right hand side determined by the overlap parameter), but were tested (and if necessary, discarded) using the $\calL$-Lebesgue stability indicator. In this work, we do not use the overlap parameter, but instead directly compute RBF-FD weights for every $\vy \in P_1$, and assess their suitability using the $\calL$-Lebesgue function at that point. This process is repeated on every stencil.

\subsubsection{An oscillation indicator}
\label{sec:natspnorm}
All RBFs have an associated reproducing kernel Hilbert space called the native space~\cite{Fasshauer:2007}. The \emph{native space semi-norm} of an RBF interpolant formed from PHS RBFs is a measure of how much the interpolant oscillates~\cite{iske2002}. In the standard interpolation setting (rather than the RBF-FD setting), PHS RBF interpolants with RBF interpolation matrix $A$ and RBF coefficient vector $\vc$ have a native space semi-norm of $|\vc^T A \vc|$~\cite{iske2002}. We now define an analogous oscillation indicator for the RBF-FD context on the stencil $P_1$ for the interpolant $s_1(\vx,\vy)$ as:
\begin{align}
\lf.\mathcal{S} (\vy)\rt|_{\vy \in P_1} = 
\lf|
\begin{bmatrix} W_1(\vy)^T & \lf(W^{\psi}_1(\vy)\rt)^T \end{bmatrix}
\begin{bmatrix}
A_1 & \Psi_1 \\
\Psi_1^T & O
\end{bmatrix}
\begin{bmatrix}
W_1 (\vy) \\
W^{\psi}_1(\vy)
\end{bmatrix}\rt|. \label{eq:mat_ind}
\end{align}
This quantity serves as a second stability indicator for the overlapped RBF-FD method. We use it to define a set $\mathbb{B}_{12}$ (analogous to $\mathbb{B}_{11}$):
\begin{align}
\mathbb{B}_{12} = \{\vy \in P_1 \mid \mathcal{S}(\vy)  \leq \mathcal{S}(\vy_1)\}.
\label{eq:ns_set}
\end{align}
While~\cite{ShankarJCP2017} found that the $\calL$-Lebesgue indicator is sufficient for low-order methods, we found that the above oscillation indicator was vital for stability in high-order methods. To better understand this indicator, we can multiply out the matrix and vectors in \eqref{eq:mat_ind} to obtain:
\begin{align}
\lf.\mathcal{S} (\vy)\rt|_{\vy \in P_1} &= \lf|W_1(\vy)^T B_{A_1} + \lf(W^{\psi}_1(\vy)\rt)^T B_{\Psi_1}\rt|,
\end{align}
where $B_{A_1}$ and $B_{\Psi_1}$ are evaluations of $\calL$ applied to the basis functions. Since application of $W_1(\vy)^T$ corresponds to approximating the action of $\calL$, this indicator can be thought of as measuring the magnitude of higher-order derivatives of the basis functions (RBFs and polynomials) as approximated by computed RBF-FD weights. While theoretical justification for this indicator is lacking, our experiments indicate that this indicator works well in conjunction with the $\calL$-Lebesgue indicator. 

\subsubsection{Parameter-Free Assembly}
\label{sec:new_assembly}
\begin{figure}[h!]
\centering
\includegraphics[scale=0.4]{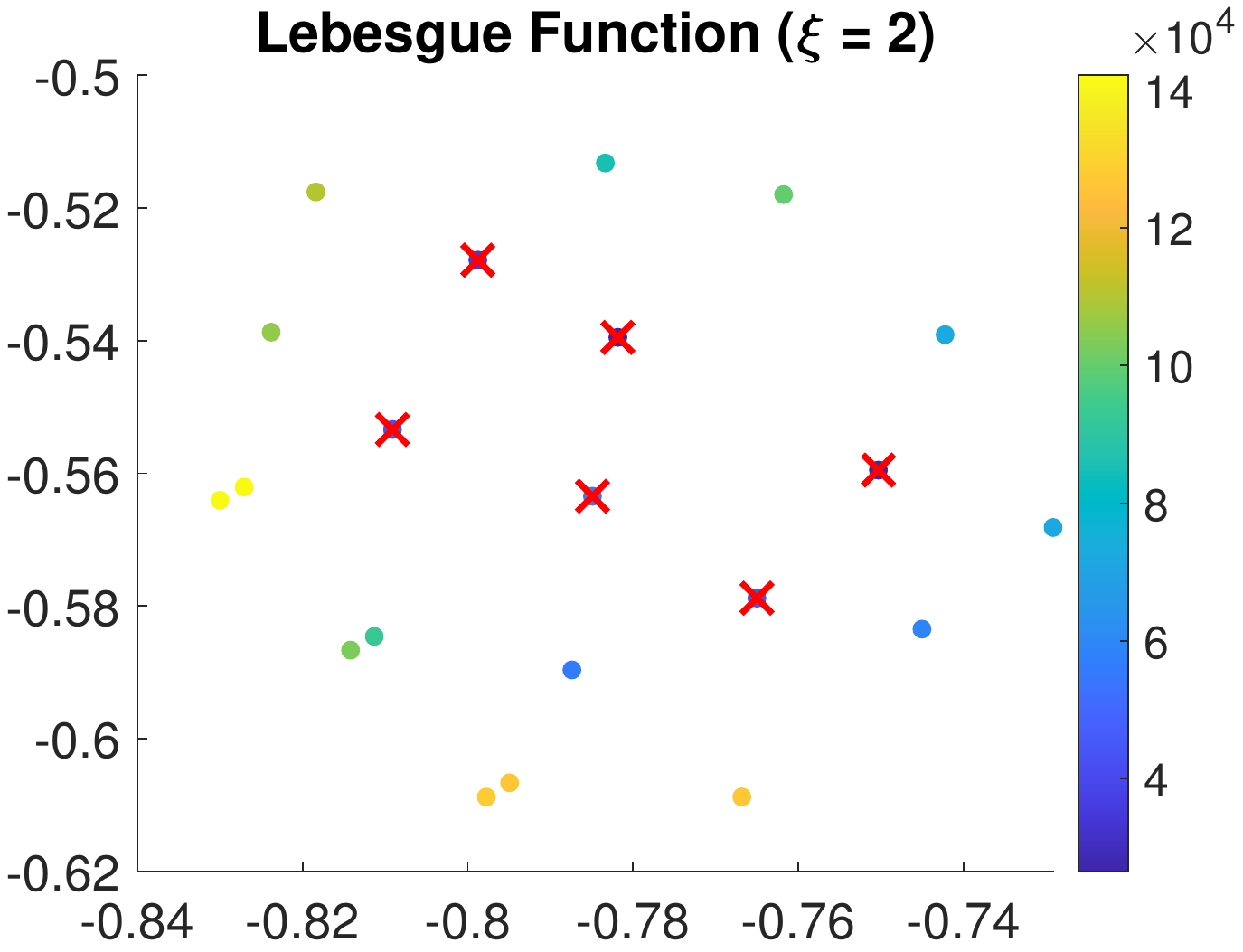}
\includegraphics[scale=0.4]{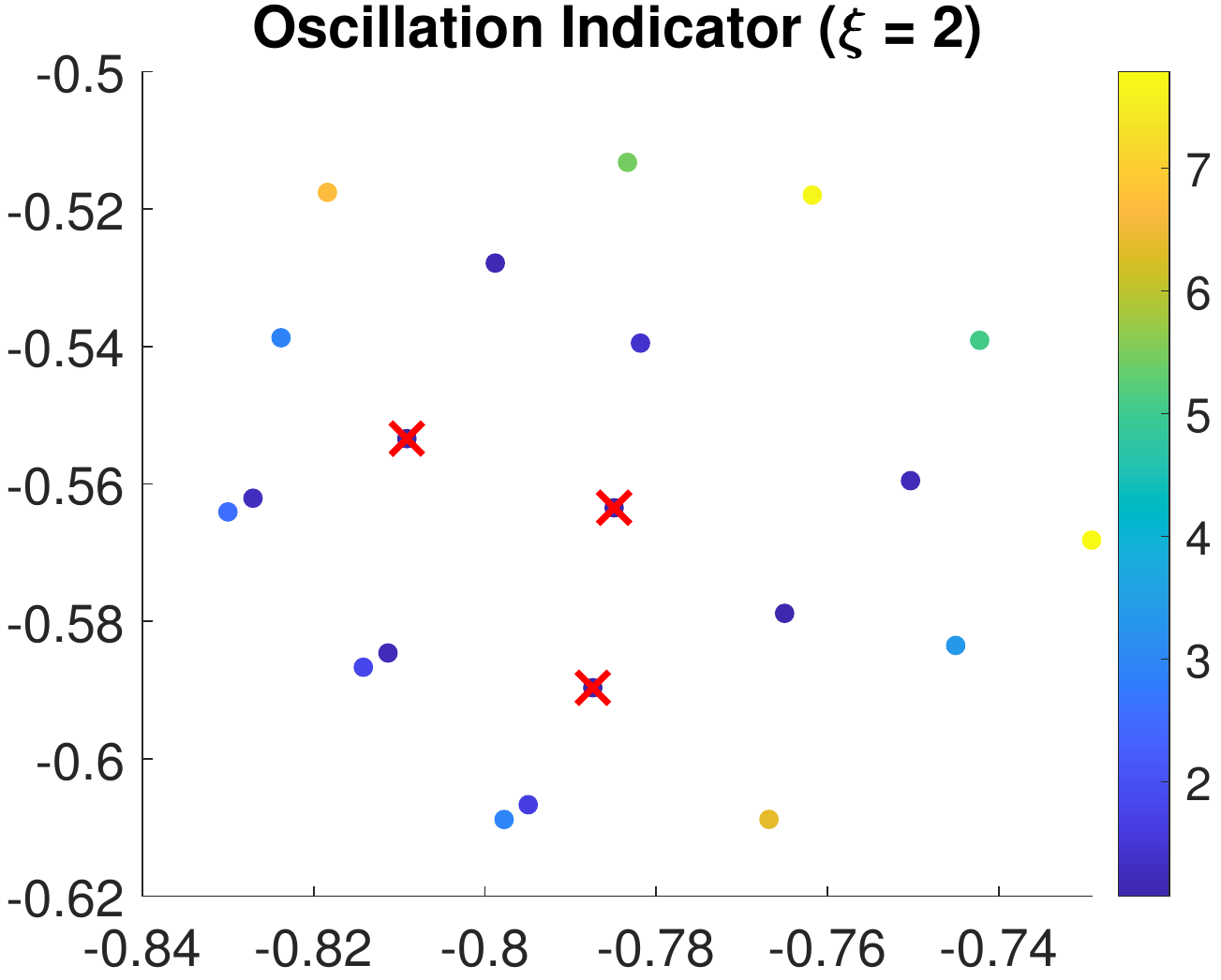}

\includegraphics[scale=0.4]{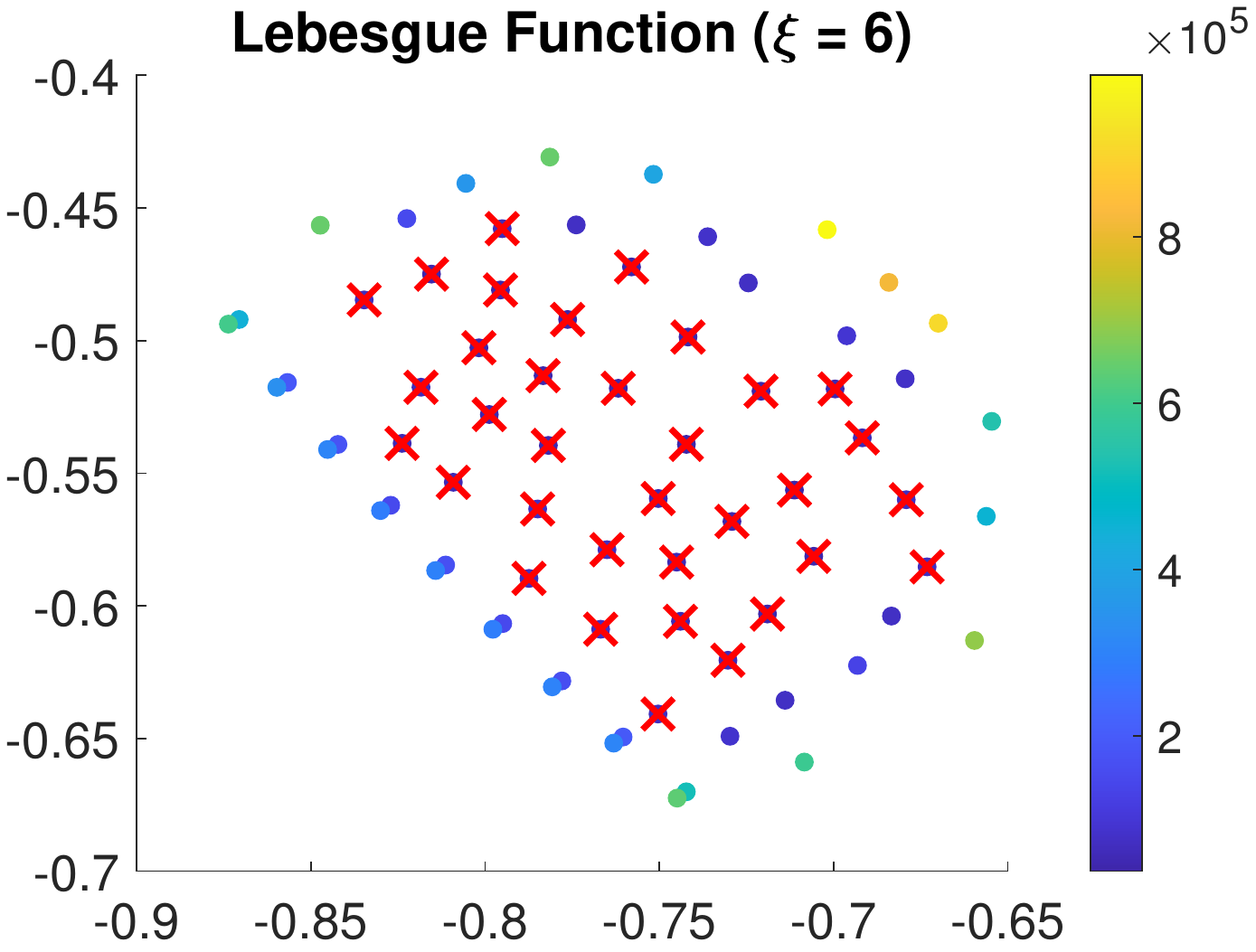}
\includegraphics[scale=0.4]{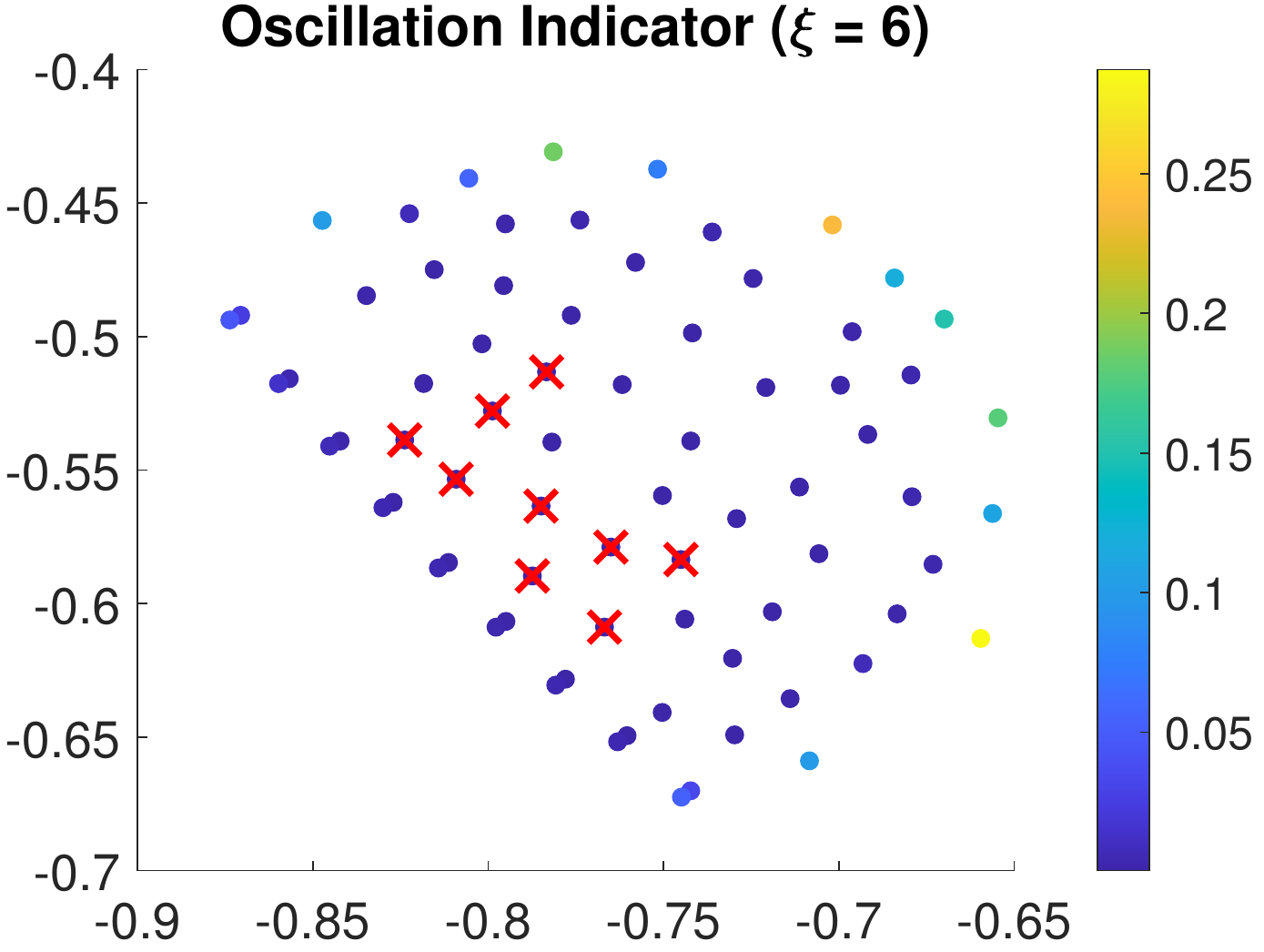}
\caption{A visualization of \eqref{eq:leb_set} (left column) and \eqref{eq:ns_set} (right column) for a stencil used to compute RBF-FD weights for the Laplacian on the unit disk ($N =4975$ nodes) with two embedded ellipses given by \eqref{eq:e1} and \eqref{eq:e2}. The colors show the function values pointwise, while the crosses indicate the points that have been selected by the indicator. The top row shows both indicators for $\xi = 2$, a second order method, and the bottom row shows the indicators for $\xi = 6$, a sixth-order method.}
\label{fig:indicators}
\end{figure}
Once both stability indicators are used to define the sets $\mathbb{B}_{11}$ and $\mathbb{B}_{12}$, we can now define a single set $\mathbb{B}_1$ of points whose weights are acceptable from the stencil $P_1$ as:
\begin{align}
\mathbb{B}_1 = \mathbb{B}_{11} \cap \mathbb{B}_{12},
\label{eq:final_seta}
\end{align}
\emph{i.e.}, the set of nodes for which we deem the RBF-FD weights as suitable is the intersection of the sets of nodes which pass both stability indicator tests. It is also useful to obtain the global indices of the points each ball $\mathbb{B}_1$:
\begin{align}
R_1 = \{\calR_1^1, \calR_2^1, \ldots, \calR_{p_1}^1\},
\label{eq:final_indset}
\end{align}
where $p_1 = |\mathbb{B}_1|$.
Note that with the definition of the sets $\mathbb{B}_k$, the stencil centers $\vy_k = \vx_{\calI^k_1}$ will automatically pass our stability tests, and an assembly algorithm based on the stencil centers is guaranteed to converge in that every point $\vx \in X$ is guaranteed to receive a set of RBF-FD weights. The indicators \eqref{eq:leb_set} and \eqref{eq:ns_set} are both visualized in Figure \ref{fig:indicators}. 

In practice, even for irregular nodes, we find that the number of stencils $N_s$ can be much smaller than the number of nodes $N$. In addition, we find that the order $\theta$ of the differential operator $\calL$ can influence the number of stencils, with higher-order operators leading to fewer stencils. This is likely due to the fact that our formula for the stencil size $n$ depends indirectly on $\theta$, and also due to the fact that larger stencils result in better behaved RBF-FD weights~\cite{ShankarJCP2017}. The full process for assembling the sparse differentiation matrix using the indicators \eqref{eq:leb_set} and \eqref{eq:ns_set} is described in Algorithm \ref{alg:new_assembly}. The new algorithm no longer requires an overlap parameter, and is therefore much more robust to irregularities in node sets. More importantly, if the node sets are changing in time (as they do in the present study), this approach obviates the need for hand-tuning.
\begin{algorithm}
\caption{Automatic Differentiation matrix assembly}
\label{alg:new_assembly}
\begin{algorithmic}[1]	
  \Statex{\bf Given}: $X = \{\vx_k\}_{k=1}^N$, the set of nodes in the domain.	
	\Statex{\bf Given}: $\calL$, the linear differential operator to be approximated.
	\Statex{\bf Given}: $n << N$, the stencil size.
	\Statex{\bf Generate}: $L$, the $N \times N$ differentiation matrix approximating $\calL$ on the set $X$.	
	\Statex{\bf Generate}: $N_s$, the number of stencils.
	\State Build a k-d tree on the set $X$ in $O(N \log N)$ operations.
	\State Initialize $g$, an array of $N$ flags set to 0.
	\State Initialize the stencil counter, $N_{s}=0$.
	\For {$k=1,N$}
		\If {$g(k)$ == 0}		
				\State Use k-d tree to determine $\{\vx_{\calI^k_1},\hdots,\vx_{\calI^k_n}\}$. Here, $\calI^k_1 = k$.
				\State Use \eqref{eq:rbf_linsys} to compute $W_k$, the $n \times n$ matrix of RBF-FD weights on the full stencil $P_1$.
				\State Find the set $\mathbb{B}_k$ using \eqref{eq:leb_set} and \eqref{eq:ns_set}.
				\State Also find the set $R_k$ by keeping track of global indices of elements of $\mathbb{B}_k$.
				\For {$i=1,n$}
					\If {$\vx_{\calR^k_i} \notin \mathbb{B}_k$}
							\State CONTINUE.
					\EndIf
					\State Set $g\lf(\calR^k_i\rt) = 1$.					
					\For{$j=1,n$}
							\State Set $L\lf(\calR^k_i,\calI^k_j\rt) = W_k(j,i)$.
					\EndFor
				\EndFor				
		\EndIf
	\EndFor
\end{algorithmic}
\end{algorithm}
}

%% file: Methods.tex
\section{A high-order meshless framework for advection-diffusion equations}

We now present a high-order semi-Lagrangian meshless method that takes advantage of the parameter-free overlapped RBF-FD formulation given in Algorithm \ref{alg:new_assembly}. Our complete method is outlined in Algorithm \ref{alg:sl_master}. Algorithm \ref{alg:sl_master} references several other algorithms and sections in this work, which will be explained later in the text. Also, since the overall framework uses a multistep method, it is important to appropriately modify the algorithm for the first two steps; we do this by using single steps of lower-order multistep methods. This section is organized as follows: first, we explain the overarching semi-Lagrangian (SL) ghost node method in Section \ref{sec:ghost}. Then, in Section \ref{sec:node_adapt}, we describe the underlying node adaption algorithm used to tackle moving embedded boundaries. Next, in Section \ref{sec:overlap_interp}, we explain how to use Algorithm \ref{alg:new_assembly} to generate overlapped local RBF interpolation stencils for use within the SL ghost node method. Parameters and error estimates are described in Section \ref{sec:error_estimates}. The preconditioner on step 16 of Algorithm \ref{alg:new_assembly} is described in Section \ref{sec:precond}. We defer a full complexity analysis of the algorithm to Section \ref{sec:comp_analysis}.
\begin{algorithm}[h!]
\caption{Semi-Lagrangian Advection-Diffusion on Moving Domains}
\label{alg:sl_master}
\begin{algorithmic}[1]	  
 \Statex{\bf Given}: $X_0$, the initial node set on the interior and boundary of a time-invariant reference domain $\Omega_0$.
	\Statex{\bf Given}: $(X_e)_0 = X_0 \cup (X_g)_0$, the initial extended node set on $\Omega_0$ containing interior, boundary, and ghost nodes.
	\Statex{\bf Given}: Seed nodes on all embedded boundaries.
	\Statex{\bf Given}: $h$, the average separation distance between nodes.
	\Statex{\bf Given}: $\xi$, the desired order of approximation of the numerical method.
	\Statex{\bf Given}: $\nu$, the diffusion coefficient.
	\Statex{\bf Given}: $\vu$, an incompressible velocity field.
	\Statex{\bf Given}: $c_0(\vx) = c(\vx,0)$, an initial condition.	
	\Statex{\bf Given}: $\triangle t$, the time step.
	\Statex{\bf Given}: $T$, the final time.	
	\Statex{\bf Given}: $g(\vx)$, the desired boundary condition such that $\mathcal{B} c = g$ over the full domain boundary.
	\Statex{\bf Generate}: $\underline{C} = \lf.C(\vx,t)\rt|_X \approx \lf.c(\vx,t)\rt|_X$, the numerical solution to \eqref{eq:lag_form_ad}.
	\State Set polynomial degree $\ell$, PHS RBF exponent $m$, and stencil size $n$ according to Table \ref{tab:params}.	
	\State Set $n_s = \lf\lfloor\frac{T}{\triangle t}\rt\rfloor$ and adjust $\triangle t$ so that $\triangle t n_s = T$.	
	\State Reconstruct embedded boundaries from seed nodes (as outlined in Section \ref{sec:node_adapt}).
	\State Modify $(X_e)_0$ to account for embedded boundaries to obtain $(X_e)(t_0)$ (as outlined in Section \ref{sec:node_adapt}).
	\State Use Algorithm \ref{alg:new_assembly} with $\calL = \delta \circ$ and $(X)(t_0)$ to obtain (localized) interpolation operator $\mathscr{I}^0$ (as outlined in Section \ref{sec:overlap_interp}).	
	\State Use Algorithm \ref{alg:new_assembly} with $\calL = \Delta$ and $\calL = \mathcal{B}$ to obtain sparse matrices $L(t_0)$ and $B(t_0$) respectively.	
	\State Set $C^0 = c_0(\vx)$.
	\For{$k = 1,\ldots,n_s$}
	\State Set $t_{n+1} = k \triangle t$.
	\State Move seed nodes on embedded boundaries using $\vu$ and RK3 with the same time-step $\triangle t$.
	\State Reconstruct embedded boundaries from seed nodes and \emph{modify} $(X_e)_0$ to obtain $(X_e)(t_{n+1})$.
	\State \emph{Update} sparse matrices $L(t_{n})$ and $B(t_{n})$ to $L(t_{n+1})$ and $B(t_{n+1})$ using Algorithm \ref{alg:op_updates}.
	\State \emph{Update} interpolation operator $\mathscr{I}^{n}$ to $\mathscr{I}^{n+1}$ using Algorithm \ref{alg:op_updates}.
	\State Trace back $(X_e)(t_{n+1})$ to $t_n$, $t_{n-1}$, and $t_{n-2}$ using \eqref{eq:sl_ode1a}--\eqref{eq:sl_ode1b} to obtain $(X_e)_d^n$, $(X_e)_d^{n-1}$, and $(X_e)_d^{n-2}$. 
	\State Compute $(C)^n_d = \mathscr{I}^n C^n$, $(C)^{n-1}_d = \mathscr{I}^{n-1} C^{n-1}$, and $(C)^{n-2}_d = \mathscr{I}^{n-2} C^{n-2}$ using the procedure outlined in Section \ref{sec:overlap_interp}.
	\State Form the preconditioner as outlined in Section \ref{sec:precond}.
	\State Form and solve the BDF3 linear system \eqref{eq:bdf3_system} (or its BDF1 or BDF2 analogues if $k=1$ or $k=2$) to obtain $C^{n+1}$.
	\State Set $C^{n-2} = C^{n-1}$, $C^{n-1} = C^n$, and $C^n = C^{n+1}$.
	\State Set $\mathscr{I}^{n-2} =  \mathscr{I}^{n-1}$,  $\mathscr{I}^{n-1} =  \mathscr{I}^n$, and $\mathscr{I}^n = \mathscr{I}^{n+1}$.
	\EndFor
\end{algorithmic}
\end{algorithm}

\subsection{A ghost node formulation}
\label{sec:ghost}
We now discuss the ghost node formulation used in Algorithm \ref{alg:sl_master}. First, we rewrite the advection-diffusion equations in \eqref{eq:adv-diff} in \emph{Lagrangian form} as 
\begin{align}
\frac{dc}{dt} &= \nu \Delta c + f(\vx,t), \vx \in \Omega(t),
\label{eq:lag_form_ad}
\end{align}
where $\frac{d}{dt} = \frac{\partial }{\partial t} + \vu \cdot \nabla$ is the \emph{material} or \emph{Lagrangian} derivative. Once the material derivative is discretized in some suitable fashion, the above equation also requires a suitable discretization of the Laplacian $\Delta$ and the boundary condition operator $\alpha \vn \cdot \nabla + \beta$ from \eqref{eq:ad-bc} on the time-varying domain $\Omega(t)$. Our approach is to use the backward differentiation formula (BDF) scheme to discretize the above equation in time, with the Laplacian, boundary conditions, and forcing terms treated implicitly. For a stable spatial discretization in the presence of derivative boundary conditions, we use the ghost node scheme outlined in~\cite{ShankarJCP2017,SFJCP2018}. We now discuss the details of this scheme in the context of domains with time-varying embedded boundaries, though the scheme is easily adapted to domains where the outer boundary also varies in time.

For a given domain $\Omega(t)$, we define a node set $X(t) = \{\vx_k(t)\}_{k=1}^{N(t)} \subset \Omega(t)$ that discretizes $\Omega(t)$. This node set is explicitly divided into a set of \emph{interior nodes} $X_i(t)$ and a set of \emph{boundary nodes} $X_b(t)$ with cardinality $N_i(t)$ and $N_b(t)$, respectively. In addition, we tile the boundary nodes some small distance in the outward normal direction to obtain a set of \emph{ghost nodes} $X_g(t)$, also of cardinality $N_b(t)$. This forms an extended node set $X_e(t) = X(t) \cup X_g(t)$, which is used extensively in Algorithm \ref{alg:sl_master}.

Our ghost node scheme involves enforcing the PDE up to and including the domain boundary $\partial \Omega(t)$, and enforcing boundary conditions at the boundary. We discretize the Laplacian using Algorithm \ref{alg:new_assembly} with $\calL \equiv \Delta$ on the extended node set $X_e(t)$. The time-varying discrete Laplacian $L(t)$ can be written in block form as:
\begin{align}
L(t) = \begin{bmatrix}
L_{ii}(t) & L_{ib}(t) & L_{ig}(t) \\
L_{bi}(t) & L_{bb}(t) & L_{bg}(t) \\
\end{bmatrix},
\end{align}
where the subscripts indicate partitions of the Laplacian corresponding to interior ($i$), boundary ($b$), and ghost ($g$) points. Notice that $L(t)$ is computed only at the interior and boundary points, but uses stencils that involve ghost points, giving the matrix dimensions of $N(t) \times (N(t) + N_b(t))$, where $N(t) = N_i(t) + N_b(t)$. We also discretize the boundary condition operator using Algorithm \ref{alg:new_assembly} with $\calL = \alpha(\vx,t) \vn \cdot \nabla + \beta(\vx,t)$. The resulting discrete boundary operator $B(t)$ can be written as:
\begin{align}
B(t) = \begin{bmatrix}
B_{bi}(t) & B_{bb}(t) & B_{bg}(t)
\end{bmatrix}.
\end{align}
This matrix has dimensions $N_b(t) \times (N(t) + N_b(t))$. The sparse matrices $L(t)$ and $B(t)$ can now be used to discretize the advection-diffusion equation. 

Let $C(\vx,t) \approx c(\vx,t)$ be the numerical solution to the advection-diffusion equation.  First, partition $\lf.C(\vx,t)\rt|_{X_e}$ into $C_i = \lf.C(\vx,t)\rt|_{X_i}$, $C_b = \lf.C(\vx,t)\rt|_{X_b}$, and $C_g = \lf.C(\vx,t)\rt|_{X_g}$. We can use these partitions to write the discretized advection-diffusion equation as:
\begin{align}
\frac{d C_i}{dt} &= \nu \lf(L_{ii}(t) C_i + L_{ib}(t) C_b + L_{ig}(t) C_g\rt) + f_i(t),\label{eq:de1}\\
\frac{d C_b}{dt} &= \nu \lf(L_{bi}(t) C_i + L_{bb}(t) C_b + L_{bg}(t) C_g\rt) + f_b(t),\label{eq:de2}\\
B_{bi}(t) C_i &+ B_{bb}(t) C_b + B_{bg}(t) C_g = g_b(t)\label{eq:de3},
\end{align} 
where $g_b(t) = \lf.g(\vx,t)\rt|_{X_b(t)}$. The above system can only be treated as a set of ODEs if a suitable discretization for the material derivative $\frac{d}{dt}$ is used. For example, for the third-order BDF scheme (BDF3)~\cite{Ascher97}, we obtain:
\begin{align}
\underbrace{\begin{bmatrix}
I - \frac{6}{11}\nu\triangle t L_{ii}^{n+1} &-\frac{6}{11}\nu\triangle t L_{ib}^{n+1} &-\frac{6}{11}\nu\triangle t L_{ig}^{n+1} \\
-\frac{6}{11}\nu\triangle t L_{bi}^{n+1} & I - \frac{6}{11}\nu\triangle t L_{bb}^{n+1} & -\frac{6}{11}\nu\triangle t L_{bg}^{n+1} \\
B_{bi}^{n+1} & B_{bb}^{n+1} & B_{bg}^{n+1}
\end{bmatrix}}_{A(t_{n+1})}
\underbrace{\begin{bmatrix}
C_i^{n+1} \\
C_b^{n+1} \\
C_g^{n+1} \\
\end{bmatrix}}_{C^{n+1}}=
\underbrace{\begin{bmatrix}
\frac{18}{11} (C_i)^n_d - \frac{9}{11} (C_i)^{n-1}_d + \frac{2}{11} (C_i)^{n-2}_d + \triangle t \frac{6}{11} f_i^{n+1}  \\
\frac{18}{11} (C_b)^n_d - \frac{9}{11} (C_b)^{n-1}_d + \frac{2}{11} (C_b)^{n-2}_d + \triangle t \frac{6}{11} f_b^{n+1}  \\
g_b^{n+1}
\end{bmatrix}}_{r^{n+1}},
\label{eq:bdf3_system}
\end{align}
where the superscripts now indicate time levels. The subscript $d$ under a variable denotes the value of that variable at an SL departure point~\cite{Xiu2002}. We  explain this in greater detail, focusing without loss of generality on the interior points. Recall that $C_i^n = \lf.C(\vx,t_n) \rt|_{X_i(t_n)}$, where $X_i(t_n)$ is the set of interior points at time level $n$. Then, the variable $(C_i)^n_d$ can be written as:
\begin{align}
(C_i)^n_d = \lf.C(\vx,t_n) \rt|_{(X_i)_d(t_n)}, 
\end{align}
where the set of \emph{departure points} $(X_i)_d(t_n) = \{ (\vx_j)_d(t_n) \}_{j=1}^{N_i(t_n)}$ is defined by solving the following set of ODEs \emph{backward} in time:
\begin{align}
\frac{d \mathbf{p}_j}{d t} &= \vu(\mathbf{p}_j,t),\label{eq:sl_ode1a}\\
\mathbf{p}_j(t_{n+1}) &= \vx_j(t_{n+1}), j=1,\ldots,N_i(t_{n+1}).
\label{eq:sl_ode1b}
\end{align}
These ODEs can be solved using a standard numerical ODE solver; \revone{we use the third-order Runge-Kutta (RK3) method}. This process is called \emph{trajectory reconstruction}~\cite{SWJCP2018}, since it reconstructs the trajectory a (fictitious) particle would take if it arrived at the nodes $X_i(t_{n+1})$; it is also alternatively referred to as the back trace procedure. In simple terms, solving the above set of ODEs for each of the $N_i(t_{n+1})$ nodes results in the set of departure points $(X_i)_d (t_n)$. In general, however, the set of departure points $(X_i)_d (t_n)$ differs from the set of interior nodes $X_i (t_n)$ at time level $n$. Consequently, $C_i^n$ must be interpolated to $(X_i)_d(t_n)$ to obtain $(C_i)^n_d$. For convenience, let us define an abstract time-dependent interpolation operator $\mathscr{I}_i(t,\cdot,\cdot)$ such that
\begin{align}
(C_i)^n_d = \mathscr{I}_i\lf(t_n, C_i^n,C_b^n\rt),
\label{eq:sl_interp1}
\end{align}
where $\mathscr{I}_i(t,\cdot,\cdot)$ is an interpolation operator that lets us interpolate fields from $X_i(t_n)$ and $X_b(t_n)$ to $(X_i)_d(t_n)$. We defer discussion of this operator to Section \ref{sec:overlap_interp}. To simplify the notation, we set 
\begin{align*}
\mathscr{I}_i\lf(t_n, C_i^n,C_b^n\rt) = \mathscr{I}_i^n \lf(C_i^n,C_b^n\rt).
\end{align*}
Note that Algorithm \ref{alg:sl_master} also requires us to compute the quantities $(C_i)^{n-1}_d$ and $(C_i)^{n-2}_d$. In analogy with \eqref{eq:sl_interp1}, these quantities can be written as:
\begin{align}
(C_i)^{n-1}_d &= \lf.C(\vx,t_{n-1}) \rt|_{(X_i)_d(t_{n-1})} = \mathscr{I}_i^{n-1}\lf(C_i^{n-1}, C_b^{n-1}\rt), \\
(C_i)^{n-2}_d &= \lf.C(\vx,t_{n-2}) \rt|_{(X_i)_d(t_{n-2})} = \mathscr{I}_i^{n-2}\lf(C_i^{n-2}, C_b^{n-2}\rt),
\end{align}
where the operators $\mathscr{I}_i^{n-1}$ and $\mathscr{I}_i^{n-2}$ interpolate quantities from $X_i$ and $X_b$ to $(X_i)_d$ at time levels $t_{n-1}$ and $t_{n-2}$ respectively. To find the departure points $(X_i)_d(t_{n-1})$ and $(X_i)_d(t_{n-2})$, we solve \eqref{eq:sl_ode1a} backward in time to levels $n-1$ and $n-2$. This requires advecting/tracing nodes $(X_i)(t_{n+1})$ backward for several steps~\cite{Xiu2002}. The same approach can be used to obtain $(C_b)_d$ values at different time levels by defining boundary interpolation operators $\mathscr{I}_b(t)$. 

Since the solution of the diffusion problem is done on the node set $X(t_{n+1})$, this approach proves far more convenient than purely Eulerian methods, which typically require spatial extrapolation to fill points that entered the domain in the current step (\emph{e.g.}, see~\cite{McCorquodaleColellaJohansen01}); this issue is completely avoided in the SL approach. \\
\begin{remark}

In Algorithm \ref{alg:sl_master}, we move the seed nodes on the embedded boundaries using an RK3 discretization of \eqref{eq:sl_ode1a} with the seed node positions in place of the nodes $\vx$, but \emph{forward} in time. In practical scenarios, it is straightforward to replace this step with updates of the type seen in the IB method~\cite{SWFKIJNMF2015}.
\end{remark}

In the following subsections, we discuss how to generate the node sets $X_e(t)$, construct and use the interpolation operators $\mathscr{I}(t)$, and efficiently \emph{update} both the discrete Laplacians $L(t)$ and the interpolation operators $\mathscr{I}(t)$.

\subsection{Node set adaptation}
\label{sec:node_adapt}
Our technique for solving the advection-diffusion equation on a time-varying domain $\Omega(t)$ involves using Algorithm \ref{alg:new_assembly} to assemble the discrete Laplacian $L(t)$ and the discrete boundary operator $B(t)$ on a \emph{time-varying node set} $X_e(t)$. In this work, we use the node generation and adaptation algorithm described in~\cite{SFKSISC2018}, adapted and optimized for moving embedded boundaries. We make this choice because the algorithm from~\cite{SFKSISC2018} is designed to \emph{locally} adapt nodes around embedded boundaries, allowing us to reuse previously computed interpolation operators and overlapped RBF-FD weights at nodes that are sufficiently far away from the moving embedded boundaries; this is discussed in greater detail in Section \ref{sec:weight_update}. We assume in the following discussion that the (irregular) outer boundary stays fixed over time, though our algorithm can be modified trivially to handle the case where the outer boundary also moves.

Let $\Omega(t)$ be a time-varying domain defined using a time-invariant reference domain $\Omega_0$ and $N_{\Gamma}$ time-varying subdomains $\{\Omega_j(t)\}_{j=1}^{N_{\Gamma}}$ so that:
\begin{align}
\Omega(t) = \Omega_0 \setminus \bigcup\limits_{j=1}^{N_{\Gamma}} \Omega_j (t).
\end{align}
The domain boundary $\Gamma(t)$ can then be written as:
\begin{align}
\Gamma(t) = \bigcup\limits_{j=0}^{N_{\Gamma}}\Gamma_j (t),
\end{align}
where $\Gamma_j(t),j=1,\ldots,N_{\Gamma}$ are the boundaries of the time-varying subdomains $\Omega_j(t)$, and $\Gamma_0 (t) = \Gamma_0$ is fixed for all time. In this setting, our goal is now to generate the set $X_e(t)$, which involves generating interior nodes $X_i(t)$, boundary nodes $X_b(t)$, and ghost nodes $X_g(t)$. Our approach is to generate node sets on $\Omega_0$ and $\Gamma_0$ prior to starting to advance the time-dependent solution, and then adapt these to account for the time-varying subdomains $\Omega_1,\ldots, \Omega_{N_{\Gamma}}$ and their boundaries $\Gamma_1,\ldots,\Gamma_{N_{\Gamma}}$.\footnote{If instead a node set is given directly on $\Omega(t_0)$, it is sufficient for our methods to precede the node adaptation procedure with a step that simply fills $\Omega_j(t_0)$, $j=1,\ldots,N_{\Gamma}$ with nodes at the initial time $t=t_0$.} The approach is as follows:
\begin{enumerate}
\item Starting from a small set of seed nodes on $\Gamma_0$, generate the node set $X_0 = (X_i)_0 \cup (X_b)_0$ for $\Omega_0 \cup \Gamma_0$ and a geometric representation for $\Gamma_0$ using Algorithm 1 from~\cite{SFKSISC2018}. Use this geometric representation to generate a set of outward unit normal vectors $\mathcal{N}_0 = \{(\vn_0)_j\}_{j=1}^{N_b}$ and use them, in turn, to generate the ghost nodes $(X_g)_0$ for the domain $\Omega_0$. This gives us the extended node set $(X_e)_0$ on the time-invariant domain $\Omega_0 \cup \Gamma_0$. 

\item Next, at any time $t$, use Algorithm 5 from~\cite{SFKSISC2018} to adapt the node set by ``turning-off'' any nodes contained in $\bigcup\limits_{j=1}^{N_{\Gamma}} \Omega_j (t)$. This requires forming a geometric representation of the boundaries $\bigcup\limits_{j=1}^{N_{\Gamma}} \Gamma_j (t)$ \revone{from a set of initally quasi-uniform seed nodes spaced $h_d$ apart}. We use the updated positions of the seed nodes and the geometric modeling technique presented in~\cite{SFKSISC2018} to form geometric representations of the embedded boundaries. This node adaptation algorithm uses normal vectors on each subdomain boundary to test whether the nodes in $(X_i)_0$ are inside or outside the domain $\Omega(t) = \Omega_0 \setminus \bigcup\limits_{j=1}^{N_{\Gamma}} \Omega_j (t)$.  \revone{The geometric modeling error is $O(h_d^8)$~\cite{SFKSISC2018}}.

\item We now need to ensure that the boundary nodes $X_b(t)$ respect the average node spacing $h$. While the outer boundary nodes are assumed to be fixed in time, the inner boundaries $\Gamma_j$ may deform over the course of a simulation. In Algorithm \ref{alg:sl_master}, we describe moving the boundaries by moving the seed nodes; this amounts to a Lagrangian description of the embedded boundaries, not unlike in the IB method~\cite{PESKIN:2002:IBM}. To avoid boundary nodes being further apart than $h$ due to this movement, \revone{we first reconstruct the embedded boundaries from their seed nodes using the previously-mentioned geometric modeling technique. We then use Algorithm 2 from~\cite{SFKSISC2018} to sample the embedded boundaries in such a way that the resulting node sets are quasi-uniform and neighboring points a distance of approximately $h$ apart}. Thus, in Algorithm \ref{alg:sl_master} (as used in this article), the seed nodes on the embedded boundaries are specified once at the beginning of the simulation (or at the introduction of an embedded boundary). In contrast, the boundary nodes are \textbf{regenerated} every time step and may vary in number. \revone{If the seed nodes drift apart, one could use the boundary representation at a given step to regenerate the seed nodes via quasi-uniform sampling with a spacing of $h_d$}. Once the the boundary nodes are obtained, the normal vectors at the boundary nodes $X_b(t)$ can be used to obtain the set of ghost nodes $X_g(t)$ (by extension in the normal direction). This generates the fully extended set $X_e(t)$. 
\end{enumerate}
\subsection{Overlapped local RBF interpolation}
\label{sec:overlap_interp}
During every step of the method, the departure points $(X_i)_d(t_{n})$, $(X_i)_d(t_{n-1})$, and $(X_i)_d(t_{n-2})$ (and their boundary counterparts) are calculated by solving \eqref{eq:sl_ode1a}. Once these departure points are calculated, the solution is interpolated to these departure points using the interpolation operators $\mathscr{I}^n$, $\mathscr{I}^{n-1}$, and $\mathscr{I}^{n-2}$. In this article, we use an analogue of RBF-FD (local RBF interpolation) to compute these interpolation operators, making them completely localized and inexpensive to compute and evaluate. In previous work on a static domain (the sphere), the authors used a separate stencil for each of the points in the domain~\cite{SWJCP2018}. \revone{Our approach here is similar, except that we now use overlapped RBF-FD to generate potentially fewer stencils.}

\revone{At a high level, our approach is to decompose the points $X(t) \subset \Omega(t)$ into stencils first, use the back trace procedure to determine which stencil a departure point lies on, then interpolate the solution from the node set $X(t)$ to the set of departure points $X_d(t) = (X_i)_d(t) \cup (X_b)_d(t)$. To automatically determine stencils, we use Algorithm \ref{alg:new_assembly} with the point evaluation operator $\delta \circ$ in place of $\calL$.  However, unlike in the case of the matrices $L(t)$ and $B(t)$, we do not assemble the resulting weights into a sparse matrix. Instead, we only store the $LU$ decompositions of the interpolation matrices. Once the stencils are determined and the $LU$ decompositions are stored, we build an acceleration structure (a k-d tree) on the stencil centers. 

For each time-step, we then determine which stencil each departure point $\lf(\vx_j\rt)_d$, $j=1,\ldots,N(t)$ lies on by finding the closest stencil center to that departure point. Suppose that some number of departure points $(\vx_j)_d$, $j=1,\ldots,N_k$ turn out to be associated with the stencil $P_k$, and for simplicity, ignore time levels for the moment. Our goal is to find $C\lf((\vx_j)_d,t\rt), j=1,\ldots,N_k$ given the values of $C(\vx,t)$ on the stencil ${P_k}$ and the $LU$ decomposition of the interpolation matrix on $P_k$. The procedure to do so is as follows:
\begin{enumerate}
\item Using the stored LU decomposition of that stencil interpolation matrix from \eqref{eq:rbf_linsys} and $\lf.C(\vx,t)\rt|_{P_k}$ as the right hand side, solve a local linear system for $O(n^2)$ operations to determine a set of RBF and polynomial coefficients $\vc^k$ on the stencil $P_k$.

\item Using these coefficients, evaluate the local interpolant on the stencil $P_k$ at each departure point $(\vx_j)_d, j=1,\ldots,N_k$ on $P_k$ to get $C((\vx_j)_d,t)$.
\end{enumerate}
This procedure is repeated for every stencil}.

Finally, it is important to note that we do not use the extended node set $X_e(t)$ to compute the interpolation stencils using Algorithm \ref{alg:new_assembly}. This would produce interpolation stencils that use unphysical values at ghost nodes to interpolate data to the departure points. Instead, the interpolation stencils are only calculated on the set $X(t) = X_i(t) \cup X_b(t)$, which does not contain ghost nodes. While this results in one-sided stencils at the domain boundary, we found this to be more stable over a wide range of Peclet numbers and boundary conditions than an approach that allowed the incorporation of ghost nodes into SL interpolation stencils. 

\subsection{Selective updates to RBF-FD weights and stencils}
\label{sec:weight_update}
\begin{algorithm}[h!]
\caption{Efficient differentiation matrix update}
\label{alg:op_updates}
\begin{algorithmic}[1]	
  \Statex{\bf Given}: $X_e(t_n)$, the $N_e(t) \times d$ matrix of nodes at time $t_n$.
	\Statex{\bf Given}: $X_e(t_{n+1})$, the $N_e(t_{n+1}) \times d$ matrix of nodes at time $t_n$.
	\Statex{\bf Given}: $\mathscr{K}_e$, a kd-tree built on $X_e(t_{n+1})$.
	\Statex{\bf Given}: $\calL$, the linear differential operator to be approximated.
	\Statex{\bf Given}: $n << N$, the stencil size.
	\Statex{\bf Given}: $L(t_n)$, the $N_e(t_n) \times N_e(t_n)$ differentiation matrix corresponding to $\calL$ at time $t_n$.
	\Statex{\bf Given}: $X_c(t_n)$, the $N_c \times d$ matrix of \emph{stencil centers} at time $t_n$.
	\Statex{\bf Given}: $N_c(t_n) = |X_c(t_n)|$, the number of stencil centers at time $t_n$.	
	\Statex{\bf Given}: $\mathscr{N}_c$, an $N_c(t_n) \times n$ matrix of nearest neighbor indices for $X_c(t_n)$ in $X_e(t_n)$.
	\Statex{\bf Given}: $\mathscr{W}_c$, a matrix mapping each row of $X_c(t_n)$ to the corresponding rows and columns in $L(t_n)$.
	\Statex{\bf Generate}: $L(t_{n+1})$, the $N_e(t_{n+1}) \times N_e(t_{n+1})$ differentiation matrix approximating $\calL$ at time $t_{n+1}$.	
	\State Build a kd-tree $\mathscr{K}_c$ on the set $X_c(t_n)$ in $O(N_c \log N_c)$ operations.
	\State Initialize $X_a = X_e(t_{n+1})$, an $N_e(t_n) \times d$ matrix of \emph{active points} (points for which weights must be computed).	
	\For {$k=1,N_c$}
		\State Use $\mathscr{K}_e$ to determine if $X_c(k,:) \in X_e(t_{n+1})$. If not, CONTINUE.		
		\State Use $\mathscr{K}_e$ to determine $n$ nearest neighbors of $X_c(k,:)$ in $X_e(t_{n+1})$. 
		\State Query $\mathscr{N}_c(k,:)$ to obtain $n$ nearest neighbors of $X_c(k,:)$ at time $t_n$.
		\If {Neighbors at $t_{n+1}$ do not match neighbors at $t_{n}$}
			\State CONTINUE.
		\EndIf
		\State Consult $\mathscr{W}_c(k)$ to obtain $p_k$ row indices (stored in $\calR_k(t_{n})$) and $n$ column indices into $L(t_n)$ (stored in $\calI_k(t_{n})$) corresponding to $X_c(k,:)$.
		\State Using $\mathscr{K}_e$, get new $p_k$ row indices (stored in $\calR_k(t_{n+1})$) and $n$ column indices (stored in $\calI_k(t_{n+1})$) into $L(t_{n+1})$.
		\State Set $L_{\calR_k(t_{n+1}),\calI_k(t_{n+1})}(t_{n+1}) = L_{\calR_k(t_{n}),\calI_k(t_{n})}(t_{n})$.
		\State Remove $(X_e)_{\calR_k(t_{n+1})}(t_{n+1})$ from the active point matrix $X_a$.
	\EndFor
	\State For all remaining points in $X_a$, use Algorithm \ref{alg:new_assembly} to compute the remaining rows of $L(t_{n+1})$.
\end{algorithmic}
\end{algorithm}
We now discuss our technique for efficiently computing the matrices $L(t)$ and $B(t)$ (in \eqref{eq:de1}--\eqref{eq:de3}), and the operators $\mathscr{I}(t)$ on the time-varying node set $X_e(t)$. \revone{Since the node sets vary in time, the number of rows and columns and also the entries of $L(t)$ and $B(t)$ are time-varying. It is instructive to enumerate all the scenarios in which these changes occur.} For the following, assume that $X_c(t)$ is the set of stencil centers, where a ``stencil center'' is always a node from the set $X(t)$.\footnote{In the standard RBF-FD method, every node from the set $X(t)$ is a stencil center. Since overlapped RBF-FD is used, there are far fewer stencil centers than nodes.} \revone{Note that any of the following scenarios could require recomputation of weights for multiple rows of the matrix $L(t)$ since overlapped RBF-FD uses the same stencil to compute weights for multiple rows at once.}
\begin{enumerate}
\item \revone{A node that was a stencil center at time $t_n$ may leave the domain at time $t_{n+1}$ by becoming covered by one of the embedded domains $\Omega_j(t)$. The set of weights associated with that stencil is no longer valid.

\item A node that was a nearest neighbor to a stencil center at time $t_n$ may no longer be a nearest neighbor to that same stencil center at time $t_{n+1}$. This could occur either if the node left the domain or another node (such as a boundary point) ended up closer. Again, the corresponding overlapped RBF-FD weights need to be recomputed.

\item The global indices of each stencil's nodes into the set $X_e(t)$ could change since the total number of nodes $N(t) + N_b(t)$ is a function of time. In this case, new global indices must be found so that the previously computed weights are copied in the appropriate rows and columns of $L(t)$.

\item New nodes may be introduced to the domain. For instance, nodes on the moving boundaries $\Gamma(t)$ are always changing, as are their ghost nodes. All these new nodes create rows in $L(t)$. Overlapped RBF-FD weights must be freshly computed for each of these nodes.}

\item In the case of $B(t)$, the differentiation matrix that enforces boundary conditions, we must therefore not only check the above cases, but also check if $\alpha$ and/or $\beta$ have changed. If they have, we need to recompute the RBF-FD weights for the boundary.
\end{enumerate}
All of the above statements for $L(t)$ also apply to $\mathscr{I}(t)$ with a minor caveat.\footnote{The only caveat is that no large sparse matrix is maintained for $\mathscr{I}(t)$, and only the stencil information and decomposed local interpolation matrices are needed.} \revone{Algorithm \ref{alg:op_updates} shows this procedure for the Laplacian, but can be trivially adapted to update the boundary condition matrix and interpolation operators also. All the nodes for which weights are not copied are marked as such, and Algorithm \ref{alg:new_assembly} is then applied to compute the weights for these nodes.}
\begin{remark}
While Algorithm \ref{alg:op_updates} describes the update procedures in terms of matrices, it is again easily adapted to a matrix-free approach where the matrices $L(t)$ and $B(t)$ are never formed, but are instead simply \emph{applied} to approximate solution vectors.
\end{remark}

\subsection{Error Estimates}
\label{sec:error_estimates}
We now discuss the error estimates for the discrete PDE \eqref{eq:bdf3_system}. To the best of our knowledge, error estimates for SL methods are only available for pure linear advection problems~\cite{FalconeFerretti}. A formal analysis in the context of a multistep method is beyond the scope of this article. However, it is possible to heuristically account for the different sources of error in \eqref{eq:bdf3_system} and combine them to get an estimate of the total error.

First, we \revone{account} for all spatial errors due to RBF interpolation or RBF-FD approximations. Recall that $h$ is a measure of average node spacing. \revone{For the overlapped RBF-FD method (including overlapped local RBF interpolation), the spatial error in approximating a differential operator of order $\theta$ is $O(h^{\ell_{\rm op} + 1 - \theta})$}, where $\ell_{\rm op}$ is the degree of the polynomial used within the RBF-FD (or local RBF) formula for that operator~\cite{DavydovSchaback2018}. In our case, $\theta=2$ for the Laplacian, $\theta=1$ for derivative boundary conditions (Neumann or Robin), and $\theta = 0$ for interpolation. Let $E_{\Delta}$ be the error in approximating the Laplacian, and $E_{\nabla}$ the error in approximating the gradient. These can now be written as:
\begin{align}
E_{\Delta} = O\lf(h^{\ell_{\Delta} -1}\rt), \ E_{\nabla} = O\lf(h^{\ell_{\nabla}}\rt),
\end{align}
where $\ell_{\Delta}$ and $\ell_{\nabla}$ are the polynomial degrees used for computing the overlapped RBF-FD weights for the operators $\Delta$ and $\nabla$ respectively.

\revone{Next, we account for the errors from the SL portion of the algorithm. \eqref{eq:bdf3_system} relies on three back traces and interpolations per step. Let $\ell_I$ be the degree of the polynomial used for overlapped local RBF interpolation. Then, the total error from the SL portion of the algorithm is given by~\cite{FalconeFerretti}:
\begin{align}
E_{SL} = O\lf( \triangle t^p + \frac{h^{\ell_I + 1}}{\triangle t}\rt).
\end{align}
where $p$ is the order of the method used to perform the backtrace ($p = 3$ for RK3)}. Finally, we have an error contribution of $O\lf(\triangle t^3 \rt)$ from the BDF3 time-stepping scheme itself. The total global error for stepping the PDE on the interior and the boundary can be written as:
\begin{align}
E_{Tot} = O\lf(h^{\ell_{\Delta} -1}\rt) + O\lf(h^{\ell_{\nabla}}\rt) + O\lf(\triangle t^3\rt) + O\lf(\triangle t^p + \frac{h^{\ell_I + 1}}{\triangle t}\rt).
\end{align}
Let $\xi$ be the desired spatial approximation order of the method. By setting $\ell_{\Delta} = \xi+1$ and $\ell_{\nabla} = \xi$, \revone{\emph{i.e.}, choosing different polynomial degrees for each differential operator}, the first two terms become $h^{\xi}$. \revone{In the last term, we set $\triangle t = O(h)$, set $\ell_I = \xi$, and $p=3$.} The error estimate therefore simplifies to
\begin{align}
E_{Tot} = O\lf(h^{\xi}\rt) + O\lf(\triangle t^3\rt).
\end{align}
These parameter choices for Algorithm \ref{alg:sl_master} are summarized in Table \ref{tab:params}.\\
\begin{remark} Our choice of $\triangle t = O(h)$ was \emph{not} motivated by stability, but rather by the desire to have low temporal error.\end{remark}
\begin{remark} While we do not have a proof of stability, it is well known that the SL framework allows for time-steps that are much larger than the CFL limit, including in the RBF context~\cite{SWJCP2018}. However, in order to prevent trajectories from crossing, it is common to set $\triangle t \leq |J|^{-1}$, where $|J|$ is the minimum pointwise Jacobian of the velocity field evaluated on the collocation node set. In this work, we find that setting $\triangle t = \frac{0.3 h}{U_{max}}$ is sufficient for both stability and accuracy, where $U_{max}$ is an estimate of the spatiotemporal maximum of $\|\vu\|_2$.\end{remark}
\begin{remark} \revone{The above analysis ignores geometric modeling errors due to moving boundaries, which are $O(h_d^{8})$. These errors are significantly smaller than the other errors in our work and can be safely ignored.}\end{remark}
\begin{table}[h!]
\centering
\begin{tabular}{ccc} \toprule
		Parameter & Meaning & Value \\ \midrule
		$\ell_I$ & Polynomial degree for Interpolation & $\xi$ \\
		$\ell_{\nabla}$ & Polynomial degree for Neumann/Robin BCs & $\xi$ \\
		$\ell_{\Delta}$ & Polynomial degree for Laplacian & $\xi+1$ \\		
		$m_{op}$ & PHS degree for operator $op$ & $\ell_{op}$ if $\ell_{op}$ is odd, $\ell_{op} - 1$ if $\ell_{op}$ is even, $m_{op} = \max(m_{op},3)$. \\		
		$n_{op}$ & Stencil size for operator $op$ & $2 {\ell_{op} + d \choose d} + 1$\\	 \bottomrule		
\end{tabular}
\caption{Table of parameters based on desired approximation order $\xi$, dimension $d$, and operator $op = \Delta$, $\nabla$, or $I$ is the operator being approximated.}
\label{tab:params}
\end{table}

%% file: LinAlg.tex
\section{Iterative methods for the implicit system}
\label{sec:linalg}

We now discuss the iterative method we use for solving the system \eqref{eq:bdf3_system}. The block matrix $A(t)$ is sparse with at most $n$ non-zero entries per row (where $n$ is the stencil size), and changes every time-step. While the update schemes outlined in Section \ref{sec:weight_update} enable fast computation of the matrices $L(t)$ and $B(t)$ that make up $A(t)$, and the right-hand-side of \eqref{eq:bdf3_system}, the solution of this time-varying linear system requires both efficient solvers and preconditioners, especially for problems in 3D. 

\subsection{An efficient saddle-point preconditioner}
\label{sec:precond}
We use the Generalized Minimum Residuals (GMRES) method~\cite{Saad2003} to solve \eqref{eq:bdf3_system}. However, it is well known that this method requires a good preconditioner to achieve faster convergence to a tolerance. We use a technique outlined in~\cite{Benzi2008}, and discuss it briefly here. For what follows, it is useful to write $A(t)$ as a $2 \times 2$ block matrix of the form:
\begin{align}
A(t) = \begin{bmatrix}
A_{11}(t) & A_{12}(t) \\
A_{21}(t) & A_{22}(t) \\
\end{bmatrix},
\end{align}
where the \revone{block $A_{11}(t)$ is given by
\begin{align}
A_{11}(t_{n+1}) = \begin{bmatrix}
I - \frac{6}{11}\nu\triangle t L_{ii}^{n+1} &-\frac{6}{11}\nu\triangle t L_{ib}^{n+1} \\
-\frac{6}{11}\nu\triangle t L_{bi}^{n+1} & I - \frac{6}{11}\nu\triangle t L_{bb}^{n+1}\\
\end{bmatrix},
\end{align}
and the other blocks naturally follow}. The matrix $A(t)$ is a (generalized) saddle-point matrix, with each of its blocks being a sparse matrix, and its \emph{Schur complement} is
\begin{align}
S(t) = A_{22}(t) - A_{21}(t) A_{11}(t)^{-1} A_{12}(t).
\end{align}
$A(t)$ has an inverse iff $S(t)$ is invertible, which in turn requires that $A_{11}(t)$ is invertible~\cite{Benzi2005}. While this appears to be the case in practice for the overlapped RBF-FD method, \revone{we now use $S(t)$ to develop a diagonal preconditioner.  Consider the block diagonal matrix $\tilde{P}(t)$ given by:
\begin{align}
\tilde{P}(t) = \begin{bmatrix} A_{11}(t) & O \\
O & S(t)\end{bmatrix}.
\end{align}
As described in~\cite{Benzi2008}, the \emph{inverse} of $\tilde{P}$ can be a reasonable preconditioner for a system involving $A(t)$. However, computing the true inverse of $\tilde{P}(t)$ requires inverting both $A_{11}(t)$ and $S(t)$. To make this process efficient, we replace $A_{11}(t)$ and $S(t)$ with diagonal matrices. First, define the diagonal matrix $\tilde{A}_{11}(t) = (A_{11})_{ii}$, $i=1,\ldots, N(t)$. Next, define the approximate Schur complement $\tilde{S}(t)$ as $\tilde{S}(t) = A_{22}(t) - A_{21}(t) \tilde{A}_{11}(t)^{-1} A_{12}(t)$. The preconditioner we use is then given by the \emph{inverse} of the diagonal matrix
\begin{align}
P(t) = \begin{bmatrix}
\tilde{A}_{11} & O \\
O & \hat{S}(t)
\end{bmatrix}, \label{eq:precond}
\end{align}
where $(\hat{S}(t))_{ii} = (\tilde{S}(t))_{ii}, i=1,\ldots,N_b(t)$. The inverse of $P(t)$ can be computed efficiently while still serving as a reasonable preconditioner to $A(t)$. This efficiency is important, since $A(t)$ (and hence $P(t)$) changes size every time-step}.
\begin{figure}[h!]
\centering
\includegraphics[scale=0.5]{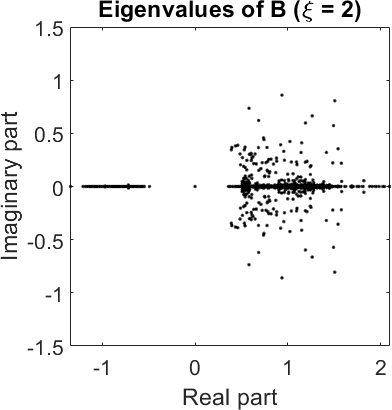} \
\includegraphics[scale=0.5]{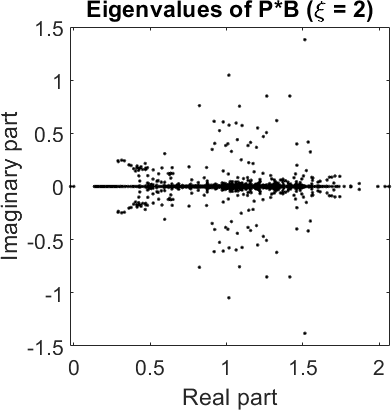}

\includegraphics[scale=0.5]{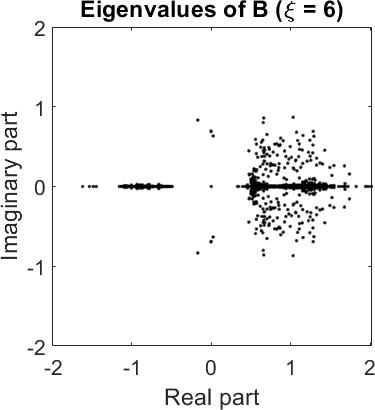} \
\includegraphics[scale=0.5]{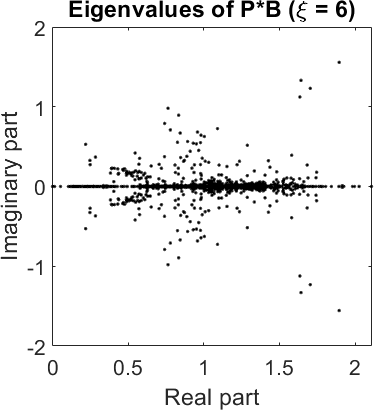}
\caption{Effect of the preconditioner $P$ from \eqref{eq:precond} on the equilibrated matrix $B$ for $\xi = 2$ (top row) and $\xi = 6$ (bottom row) for a 2D advection-diffusion problem with $\nu = 1$ at time $t=4 \times 10^{-6}$.}
\label{fig:precond}
\end{figure}
\revone{In practice, we first equilibrate $A(t)$ using Matlab's built-in equilibrate function~\cite{MatlabEquilibrate} to permute and rescale $A$ so that its off-diagonal entries are not greater than 1 in magnitude and its diagonal entries are only 1 or -1. We then form $P(t)$ using the blocks of the resulting permuted matrix $B(t)$. The resulting spectra of $B$ and $PB$ are shown in Figure \ref{fig:precond} for $\xi = 2$ and $\xi = 6$ (with $\nu = 1$). Our preconditioner has the effect of shifting most of the eigenvalues of $B$ to one side of the imaginary axis. We found in practice that this decreased the number of GMRES iterations by half. This preconditioner is clearly not optimal, and can be improved by adding more blocks from the original matrix to $P(t)$. We leave a deeper investigation of preconditioners for future work.}

\subsection{A good guess for GMRES}
\label{sec:guess}
For GMRES to converge rapidly in solving \eqref{eq:bdf3_system}, \revone{it is also important to supply a good initial guess to the solver}. In previous work, the authors have used the solution $C^n$ from the previous time level to good effect~\cite{SFJCP2018,ShankarJCP2017}. However, for a time-varying domain $\Omega(t)$, the solution $C^n$ lies on the domain $\Omega(t_n)$, while the solution $C^{n+1}$ lies on $\Omega(t_{n+1})$. Consequently, the lengths of the vectors $\lf. C(\vx,t_n)\rt|_{\vx\in X_e(t_{n})}$ and $\lf. C(\vx,t_{n+1})\rt|_{\vx\in X_e(t_{n+1})}$ are different, since the node sets $X_e(t) \subseteq \Omega(t)$ themselves vary in time. \revone{On the other hand, we see that the vectors $(C_i)^n_d$ and $(C_b)^n_d$ generated from the SL trajectory reconstruction possess cardinality $N_i(t_{n+1})$ and $N_b(t_{n+1})$, respectively. We therefore use these as initial guesses for $C_i^{n+1}$ and $C_b^{n+1}$, respectively. To obtain an initial guess for $C_g^{n+1}$, we approximate the value of $C_g$ at the ghost node departure points obtain by back tracing from the ghost nodes at $t_{n+1}$, \emph{i.e.}, we compute $(C_g)^n_d \approx \lf.C(\vx,t)\rt|_{\lf(X_g\rt)_d\lf(t_n\rt)}$. In practice, we compute this quantity by evaluating the local overlapped RBF interpolants from Section \ref{sec:overlap_interp} at the ghost departure points. We found that use of this guess vector accelerated GMRES when compared to using the zero vector as a guess (the number of iterations were decreased by an order of magnitude in 3D).}  

%% file: CompAnalysis.tex
\section{Complexity Analysis}
\label{sec:comp_analysis}

We now analyze the computational complexity of Algorithm \ref{alg:sl_master} in terms of (1) the preprocessing (steps 1-7), and (2) the actual time-stepping loop (steps 8-20). Steps 1, 2, and 7 merely involve function evaluations and will be ignored for the purposes of simplicity. For the remainder of the section, we assume without loss of the generality that the stencil size $n$ and the polynomial degree $\ell$ are the same for all operators. All our estimates nevertheless hold true in the worst case sense.

\subsection{Preprocessing complexity}
\label{sec:pre_comp}

Consider first step 3, which involves interpolation to form the geometric model of the embedded boundaries, and step 4 which involves removal of nodes from the set $(X_e)_0$. Recall that we have $N_{\Gamma}$ boundaries, and assume that there are $N_d$ seed nodes on each boundary. Since the geometric model from~\cite{SFKSISC2018} involves a dense matrix solve, the total cost of \emph{forming} this geometric model is $O(N_{\Gamma} N_d^3)$. In addition, step 4 involves both evaluation of the geometric model to generate boundary points and a subsequent modification of $(X_e)_0$. Let $N_0 = |(X_e)_0|$. Following Section 4.2 of~\cite{SFKSISC2018}, the costs of steps 3 and 4 can be rewritten as $O(N_{\Gamma} N_0)$, \emph{i.e.,} linear in $N_0$ (this relationship can be derived by rewriting the $N_d^3$ term in terms of $N_0$~\cite{SFKSISC2018}). This is applicable in both 2 and 3 dimensions.

The complexity analysis of steps 5 and 6 \revone{depends on the behavior of Algorithm \ref{alg:new_assembly}}, which in turn depends on the type of node set and the differential operator being approximated. Given the polynomial degree $\ell$, the number of polynomial terms $M = {\ell + d \choose d}$, and the stencil size $n = 2M+1$,  Algorithm \ref{alg:new_assembly} will produce a worst case complexity of $ O( (n + (n-1)/2)^3 N_0)$, where the cubic term comes from the LU decomposition of \eqref{eq:rbf_linsys} and $N_0$ is the number of points on the reference domain. However, in practice, far fewer linear systems than $N_0$ will be solved, especially as $n$ and $\ell$ are increased, since each stencil will be used to compute weights for more than one of the $N_0$ nodes. \revone{To estimate this cost, assume that for each $n$-node stencil}, Algorithm \ref{alg:new_assembly} computes RBF-FD weights for $\frac{n}{\kappa}$ nodes, where $1 \leq \kappa \leq n$. \revone{Letting} $N_s$ be the number of stencils generated by Algorithm \ref{alg:new_assembly}, the total cost of Algorithm \ref{alg:new_assembly} can be written as
\begin{align}
T_{assembly} = O\lf(N_s \lf( (n + (n-1)/2)^3 + \frac{n}{\kappa} (n+(n-1)/2)^2 \rt)\rt),
\end{align}
where the first term \revone{again} corresponds to the cost of LU decomposition, and the second term to the cost of back substitutions. We can now bound $N_s$ in terms of $N_0$ and $\kappa$. To do so, we need only realize that
\begin{align}
N_s \approx N_0 \frac{\kappa}{n}, \label{eq:Ndelta}
\end{align}
where $\kappa$ is typically closer to 1 than to $n$, making the constant $\frac{\kappa}{n}$ quite small. Thus, the assembly cost can be rewritten as
\begin{align}
T_{assembly} =  O\lf(\frac{\kappa}{n} N_0 \lf((n + (n-1)/2)^3 + \frac{n}{\kappa} (n+(n-1)/2)^2 \rt)\rt),
\end{align}
which is significantly smaller than if \revone{$\kappa = n$ as in the standard RBF-FD method}. In practice, as $n$ and $\ell$ increase, $\kappa$ decreases because larger stencils allow us to retain more weights \revone{per stencil~\cite{ShankarJCP2017}}. The only difference between steps 5 and 6 is that the values of $n$ and $\ell$ are potentially different for the Laplacian versus point evaluation.

\revone{Steps 5 and 6 use Algorithm \ref{alg:new_assembly}, which also require k-d trees to be built on the node sets for a cost of $O(N_0 \log N_0)$. The cost of} searching them for $n$ nearest neighbors on $N_s$ stencils is $O(n N_s \log N_0)$. Using \eqref{eq:Ndelta}, we can write the total preprocessing cost as
\begin{align}
T_{preprocessing} = O\lf(N_{\Gamma} N_0 + \lf((n + (n-1)/2)^3 + \frac{n}{\kappa} (n+(n-1)/2)^2 \rt) \frac{\kappa}{n} N_0 + N_0 \log N_0 + \kappa N_0 \log N_0\rt),  
\label{eq:prep_cost}
\end{align}
where again $\kappa$ is a number closer to 1 than to $n$. \revone{We have observed in practice that the second term dominates this cost and scales as $O(N_0)$ for a given $n$ and $\kappa$. To control efficiency, $\kappa$ could be explicitly introduced as an input to Algorithm \ref{alg:new_assembly}, but we leave this approach for future work.}

\subsection{Time-stepping complexity}
\label{sec:time_comp}

We now estimate the complexity for steps 8-20 which are carried out in each timestep. We ignore steps 8, 9, 18, 19, and 20, as they are trivially estimated. Instead, we focus on the cost of a single time-step. Step 10 is done using the RK3 method, which has three stages, but the cost only depends on the number of seed nodes and the total number of embedded boundaries. This cost is clearly $O(N_{\Gamma} N_d)$. Step 11 is the same as the preprocessing steps 3 and 4, and therefore has a cost of $O(N_{\Gamma} N_0)$. 

Steps 12 and 13 \revone{involve a mix of \emph{copying} old information and \emph{computing} new information}, as shown in Algorithm \ref{alg:op_updates}. Estimating the cost of these steps requires an analysis of Algorithm \ref{alg:op_updates}, combined with a slightly modified analysis for Algorithm \ref{alg:new_assembly}. Within Algorithm \ref{alg:op_updates}, letting $N_e = |X_e(t)|$, a k-d tree is first formed on $X_e(t)$ for $O(N_e \log N_e)$ operations in step 1. Steps 2 and 13 can be ignored as this can simply be done with Boolean flags. Steps 6, 10, and 12 can be done in constant time. \revone{This leaves the following steps:
\begin{itemize}
\item Step 4 can be done in $O(\log N_e)$ operations (with possible early termination).
\item Step 5 has a cost of $O(n \log N_e)$.
\item Step 7 can be done in practice by computing $\|p_{new} - p_{old}\|_2$, where $p_{new}$ are the positions of neighbors at $t_{n+1}$, and $p_{old}$ the positions at $t_n$; this costs $O(n)$ operations. Steps 8 and 9 can be ignored.
\item Step 11 involves a kd-tree search for a cost of $O(n \log N_e)$.
\end{itemize}
}Each of these steps could, in the worst case, be executed $N_s$ times (once for each stencil center). At the end of step 14, some fraction of the rows of $L(t_{n+1})$, $B(t_{n+1})$, and some of the local operators that constitute $\mathscr{I}^{n+1}$ will have been computed. Let this fraction be $\tau$, so that step 15 now only operates on $(1-\tau)N_e$ nodes. The cost of step 15 can therefore be written based on the analysis from Section \ref{sec:pre_comp} \revone{as
\begin{align}
T_{step \ 15} = O\lf(\lf((n + (n-1)/2)^3 + \frac{n}{\kappa} (n+(n-1)/2)^2 \rt) \frac{\kappa}{n} (1-\tau) N_e + \kappa(1-\tau) N_e \log N_e\rt).
\label{eq:step15}
\end{align}
The total cost of Algorithm \ref{alg:op_updates} can thus be written using the above list, \eqref{eq:Ndelta}, and \eqref{eq:step15} as
\begin{align}
T_{updates} = O\lf(\lf(1 + \frac{\kappa}{n} + 2\kappa\rt) N_e \log N_e  + \lf(3 \frac{\kappa}{n} + \kappa\rt) N_e\rt) + T_{step \ 15}.
\end{align}
In practice, $\tau \approx 1$ and $\kappa = O(1)$, allowing us to simply write
\begin{align}
T_{updates} = O(N_e \log N_e),
\end{align}
\emph{i.e.}, k-d tree lookups dominate the cost of updating the differentiation matrices $L(t)$ and $B(t)$, and the interpolation operators $\mathscr{I}(t)$. In contrast}, if $\tau$ is small, most of the rows of $L$ and $B$, and most of the local interpolants comprising $\mathscr{I}$ have to be recomputed. This concludes the analysis of steps 13 and 14 of Algorithm \ref{alg:sl_master}.

Steps 14 and 15 are the \revone{SL steps in the algorithm. Step 14 of Algorithm \ref{alg:sl_master} involves 3 back traces costing $O(N_e)$}, while step 15 \revone{implicitly involves a k-d tree lookup for each of the $N_e$ points to find the correct local interpolation stencil, one back-solve on that stencil, and at least one evaluation on that stencil}. In the worst case, we have $N_s = \frac{\kappa}{n}N_e$ stencils, exactly one back substitution per stencil for a cost of $O((n+(n-1)/2)^2)$, and exactly one evaluation per stencil also for a cost of $O((n+(n-1)/2)^2)$. The lookup costs $O\lf(N_e \log N_e\rt)$ operations. The total cost of Steps 14 and 15 of Algorithm \ref{alg:sl_master} is therefore
\begin{align}
T_{SL} = O\lf( N_e \log N_e + 2\frac{\kappa}{n}N_e (n+(n-1)/2)^2\rt),
\end{align}
which is dominated by the second term. In practice, however, it is possible for our batching technique to result in multiple evaluations per stencil and fewer back-substitutions.

Step 16 of Algorithm \ref{alg:sl_master} involves forming the preconditioner as outlined in Section \ref{sec:precond}. The cost of this step is dominated by the cost of forming the approximate Schur complement, which in turn involves sparse matrix multiplications of matrices with at most $n$ non-zero entries per row. This incurs a cost of \revone{$O(n^2 N_e)$}. Finally, the complexity of step 17 of Algorithm \ref{alg:sl_master} \revone{is non-deterministic as it} involves the convergence of the preconditioned GMRES method bootstrapped with a guess. Rather than attempt to estimate the complexity, we show results in terms of number of GMRES iterations in Section \ref{sec:results}. Thus, the total cost of each time-step of Algorithm \ref{alg:sl_master} (assuming $\tau \approx 1$) is $T_{step} = T_{updates} + T_{SL} + O\lf(n^2 N_e\rt) + {\rm Cost \ of \ GMRES}$.
\begin{align}
\implies T_{step} &= O\lf(2N_e \log N_e + 2\frac{\kappa}{n}N_e (n+(n-1)/2)^2 + n^2 N_e\rt) + {\rm Cost \ of \ GMRES}.
\end{align}
Depending on which terms dominate above, this cost is either linear or log-linear in $N_e$ (provided the cost of GMRES is kept low with a good preconditioner and initial guess). In Section \ref{sec:timings}, we demonstrate that the complexity of our method is indeed very close to linear in both the preprocessing and per-step costs, with a very large improvement in per-step costs (over the preprocessing costs) due to the reuse of weights over several steps.

%% file: Results.tex
\section{Numerical Results}
\label{sec:results}
We test our numerical framework via convergence studies on the forced advection-diffusion equation with two Peclet numbers: 1 and 1000. The forcing term is selected to maintain a prescribed solution for all time, and the prescribed solution is used to test spatial convergence rates. We solve this test problem on irregular 2D and 3D domains with boundaries moving at the fluid velocity.  Given a true solution $c(\vx,t)$ and a numerical solution $C(\vx,t)$, we compute relative $\ell_2$ errors at the final time $t=0.5$ on the node set $X$ as $e_{\ell_2} = \frac{\|c_X - C_X\|_2}{\|c_X\|_2}$. In addition to relative $\ell_2$ errors, we also report the average number of GMRES iterations per step for each value of the node count $N$ and approximation order $\xi$. Finally, we verify the complexity estimates in Section 6 via timings, and also show a comparison of computational cost and accuracy as a function of the approximation order. In all cases, we set the GMRES tolerance to $\min (0.1 h^{\xi},10^{-7})$.
\input{results2d}
\input{results3d}
\input{timings}

%% file: results2d.tex
\subsection{Forced advection-diffusion in a time-varying 2D domain}
\begin{figure}[h!]
\centering
\subfloat[]
{
	\includegraphics[scale=0.5]{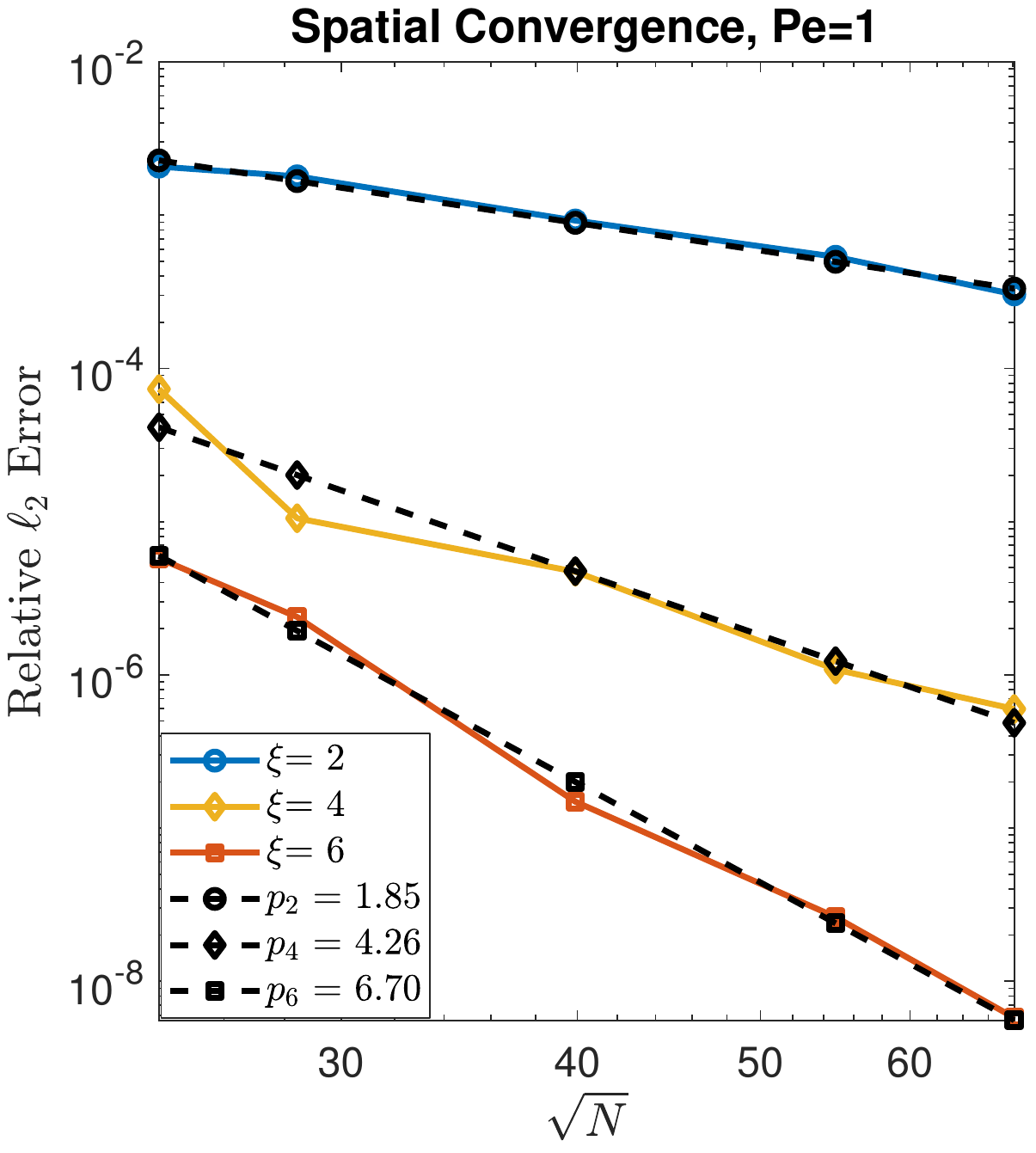}
	\label{fig:res2d_11}	
}
\subfloat[]
{
	\includegraphics[scale=0.5]{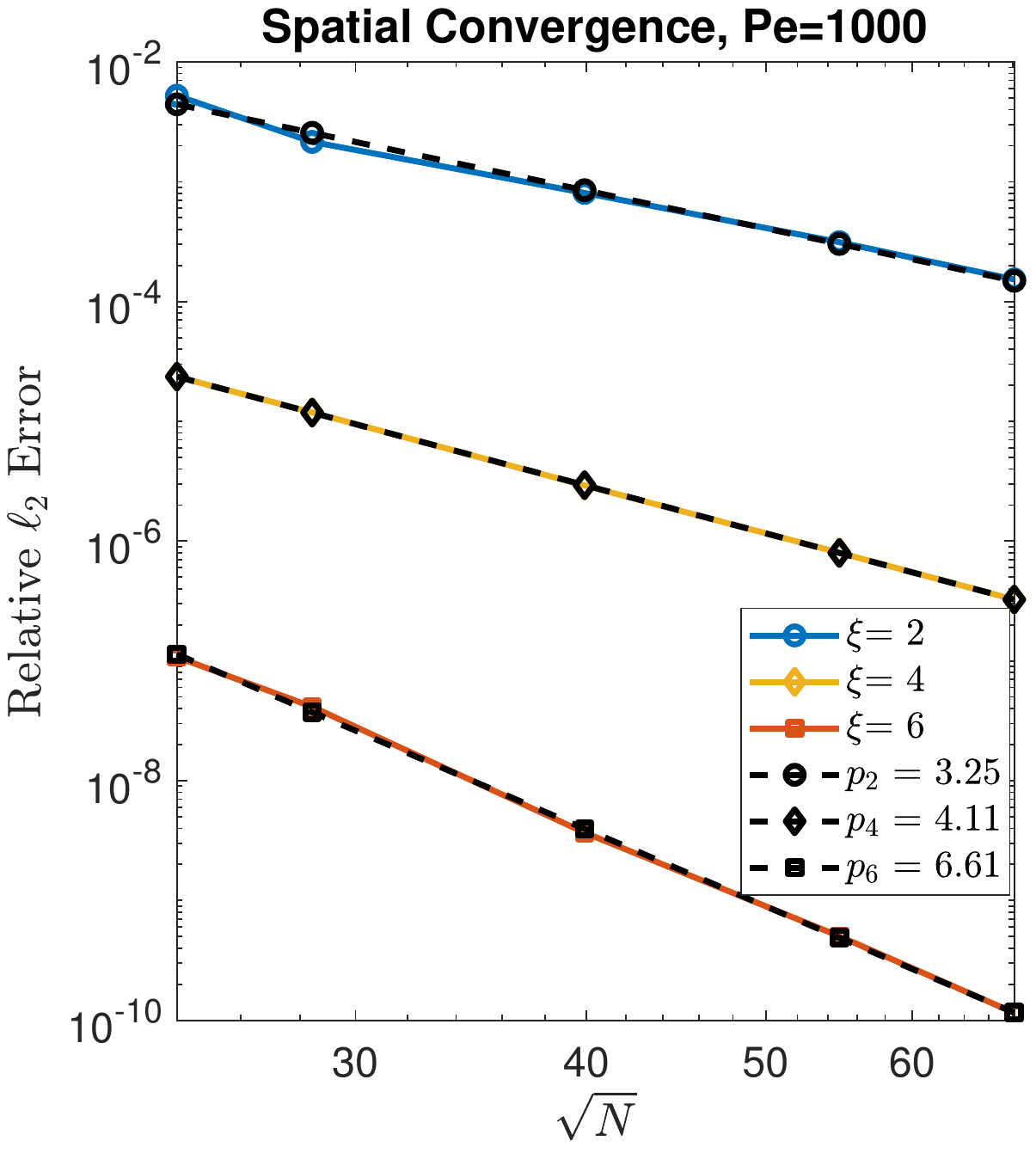}
	\label{fig:res2d_12}	
}

\subfloat[]
{
	\includegraphics[scale=0.5]{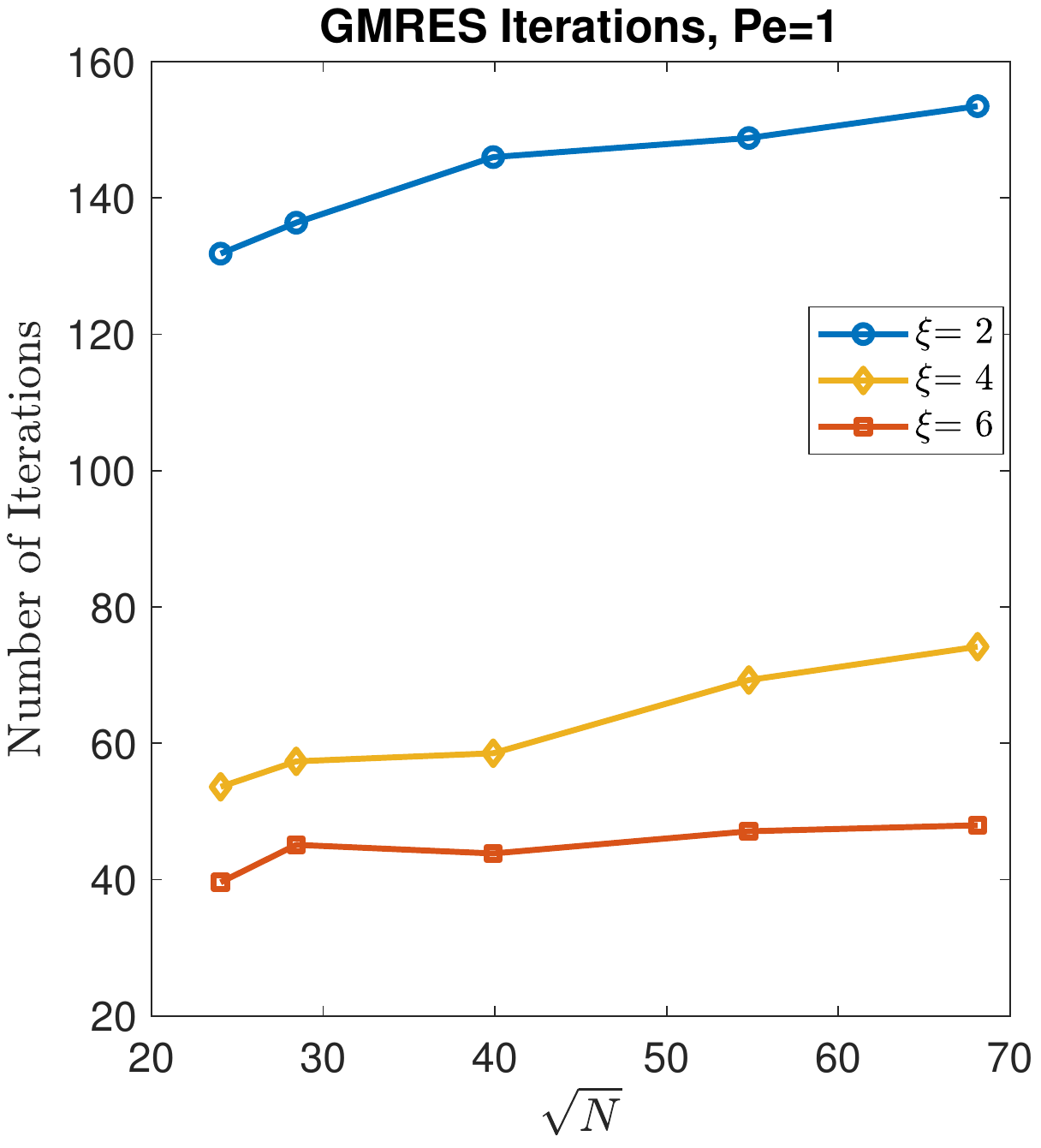}
	\label{fig:res2d_21}	
}
\subfloat[]
{
	\includegraphics[scale=0.5]{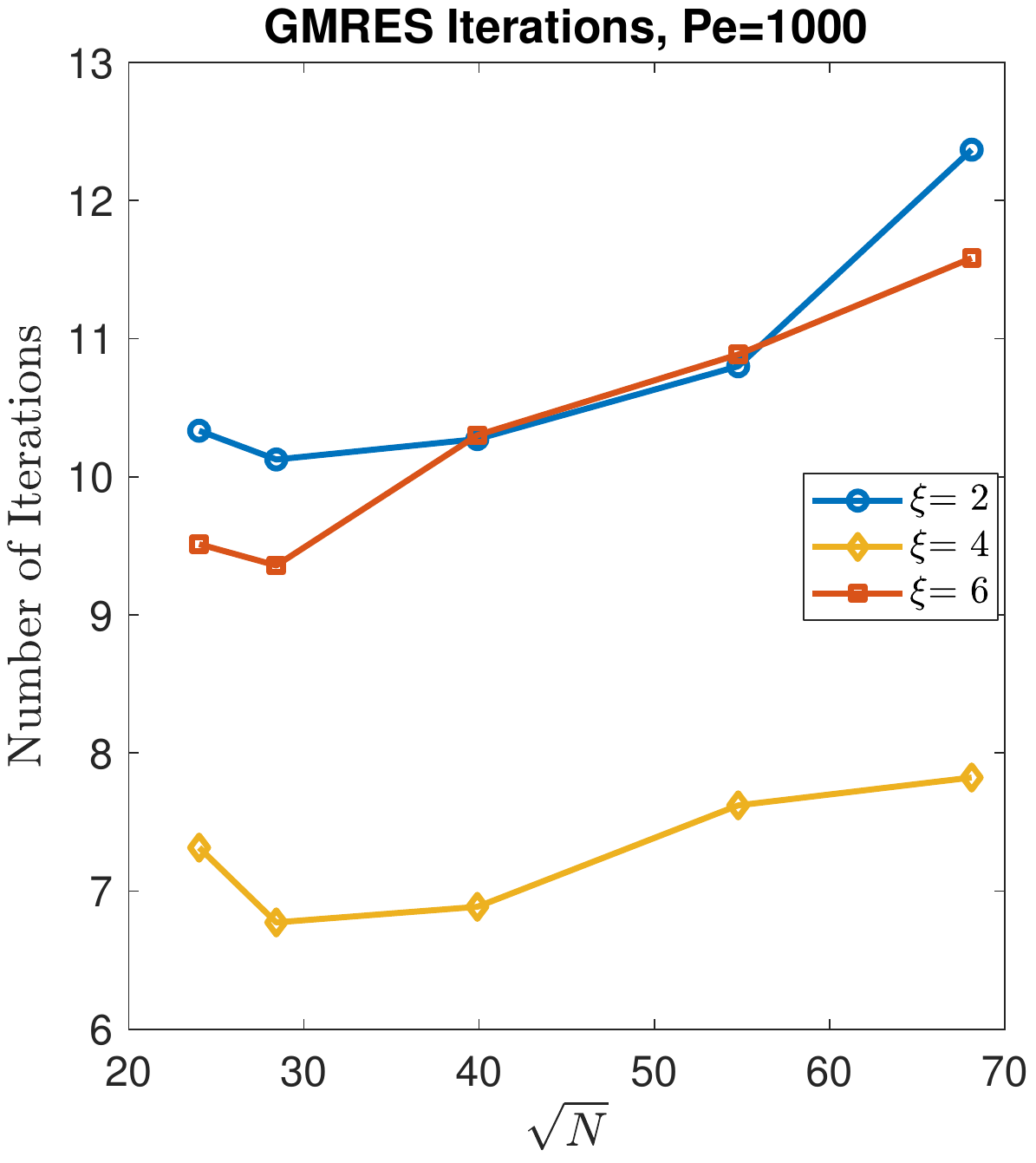}
	\label{fig:res2d_22}	
}
\caption{Top row: relative $\ell_2$ error vs $\sqrt{N}$ as a function of approximation order $\xi$ for forced advection-diffusion on the unit disk. The dashed lines are lines of best fit indicating the slope (and hence convergence rate). Bottom row: average number of GMRES iterations per time-step as a function of $\sqrt{N}$ and $\xi$.}
\label{fig:res2d}
\end{figure}
The first test involves solving the advection-diffusion equation on an irregular time-varying 2D domain. The initial domain is the unit disk with two embedded ellipses defined in parametric form as:
\begin{align}
E_1&: \ x = 0.4 \cos(\mu), y = -0.5 + 0.2 \sin(\mu),\label{eq:e1}\\
E_2&: \ x = 0.1 \cos(\mu), y = 0.2 \sin(\mu),\label{eq:e2}
\end{align}
where $-\pi \leq \mu < \pi$. The initial simulation domain is then given as $
\Omega(t_0) = \{ \mathbb{B}^2 \setminus \lf(E_1 \cup E_2\rt)\}$, \emph{i.e.}, the portion of the domain outside the ellipses but within the unit disk. We use 20 seed nodes on the boundary of each ellipse, and reconstruct the moving boundary from these seed nodes using the (parametric) periodic degree-7 polyharmonic spline RBF within the geometric model developed by the authors~\cite{SFKSISC2018}. The manufactured solution we use is given by
\begin{align}
c(\vx,t) = 1 + \sin(\pi x) \cos(\pi y) \sin(\pi t),
\end{align}
and the incompressible velocity field $\vu(\vx,t) = [u,v]$ is given by
\begin{align}
\vu(\vx,t) &= \sin\lf(\pi \|\vx\|_2^2 \rt) \sin(\pi t) [y,-x].
\end{align}
This velocity field vanishes on the boundary of the disk, and the 20 seed nodes on the embedded domain boundaries are advected with the velocity $\vu$, thereby leading to deformation of the embedded ellipses. We use a pure Neumann boundary condition operator $\mathcal{B}$ ($\alpha = -\nu$ and $\beta = 0$), and the right hand side $g(\vx,t)$ is given by applying $\mathcal{B}$ to the prescribed $c(\vx,t)$.

We measure errors in our numerical solution against the prescribed $c$. To obtain Peclet numbers $\rm{Pe} =1$ and $\rm{Pe}=1000$, we set $\nu = 1$ and $\nu = 10^{-3}$, respectively. We set the time-step as mentioned previously. We simulate the PDE to time $t=0.5$ using Algorithm \ref{alg:sl_master}.

The results for $\xi = 2,4,6$  are shown in Figures \ref{fig:res2d}(a) and \ref{fig:res2d}(b), plotted as a function of $\sqrt{N}$ (proportional to $1/h$). The results show that our predicted spatial convergence rate of $h^{\xi}$ roughly holds under refinement for both Peclet numbers $\rm{Pe} = 1$ and $\rm{Pe} = 1000$; however, the errors are much lower for $\rm{Pe} = 1000$. In addition, Figures \ref{fig:res2d}(c) and \ref{fig:res2d}(d) show the average number of GMRES iterations per timestep for each of these simulations, again as a function of $\sqrt{N}$. For $\rm{Pe}=1$, increasing $\xi$ results in fewer GMRES iterations, possibly due to the interplay of $\xi$ with the guess and preconditioner. For $\rm{Pe}=1000$, the iteration counts are more erratic as $\xi$ is increased, but it is clear that far fewer GMRES iterations are required for this value of $\rm{Pe}$, since a smaller value of $\nu$ was used (improving the conditioning of the time-stepping matrix).

%% file: results3d.tex
\subsection{Forced advection-diffusion in a time-varying 3D domain}
\label{sec:res3d}
\begin{figure}[h!]
\centering
\subfloat[]
{
	\includegraphics[scale=0.5]{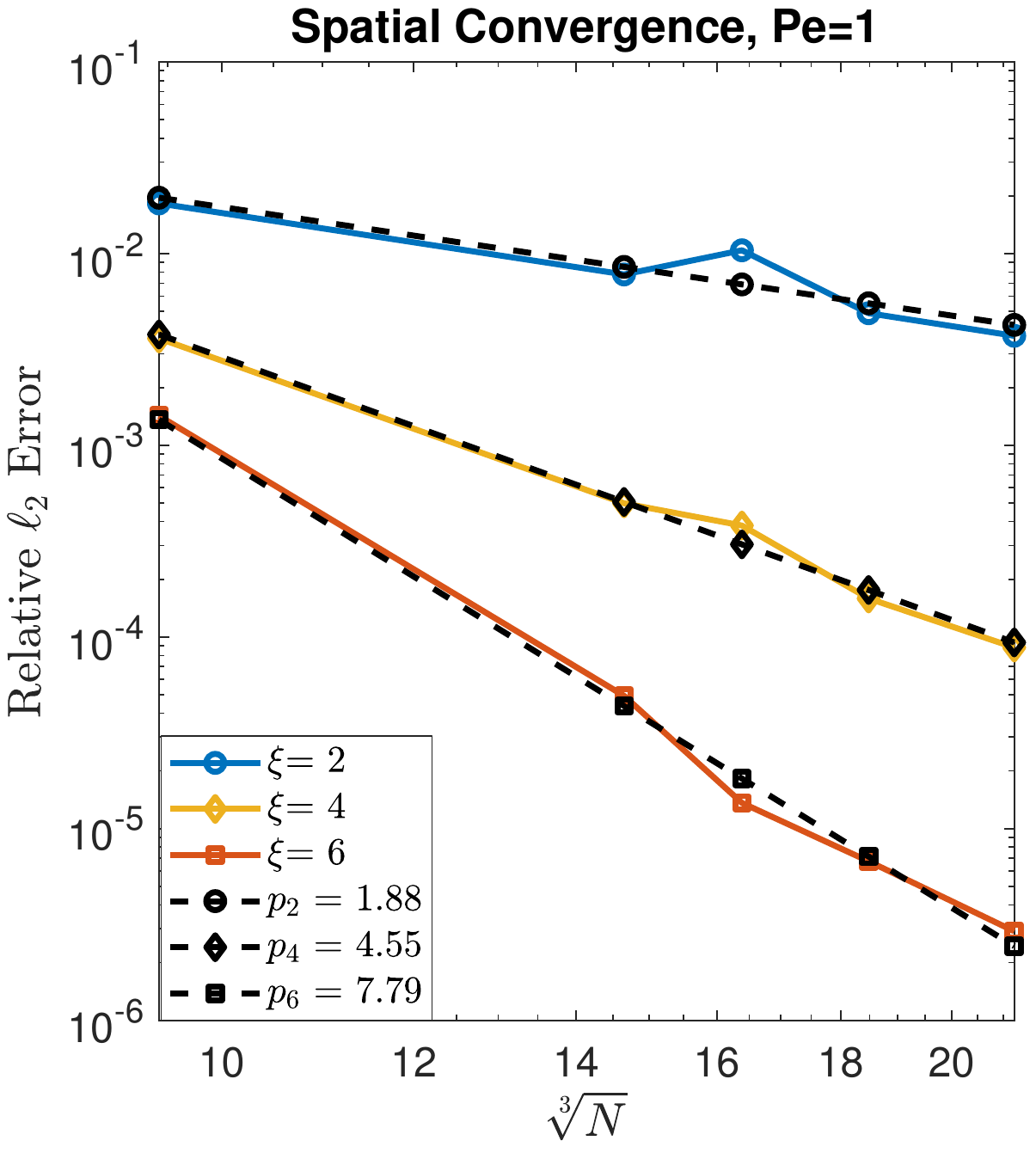}
	\label{fig:res3d_11}	
}
\subfloat[]
{
	\includegraphics[scale=0.5]{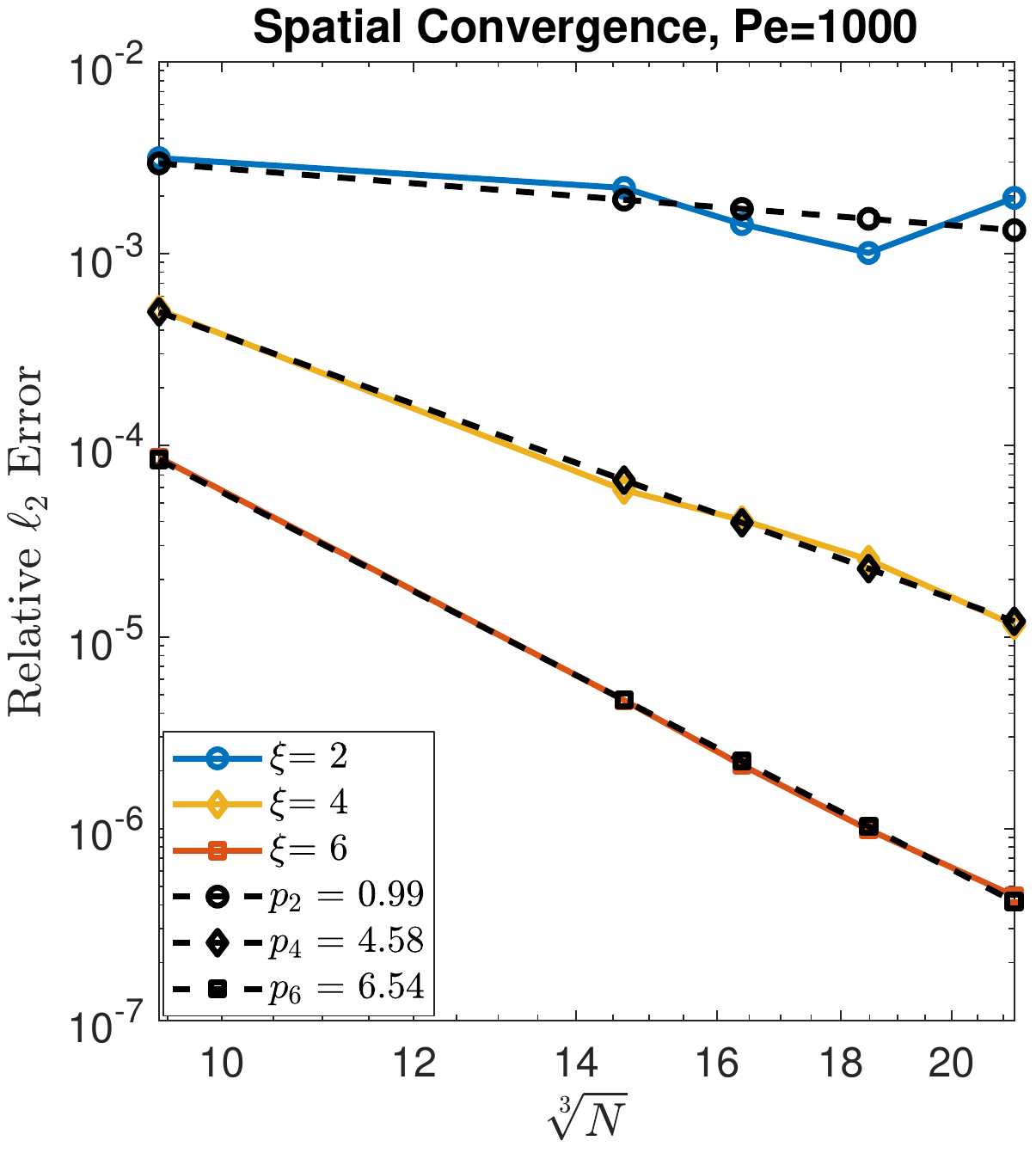}
	\label{fig:res3d_12}	
}

\subfloat[]
{
	\includegraphics[scale=0.5]{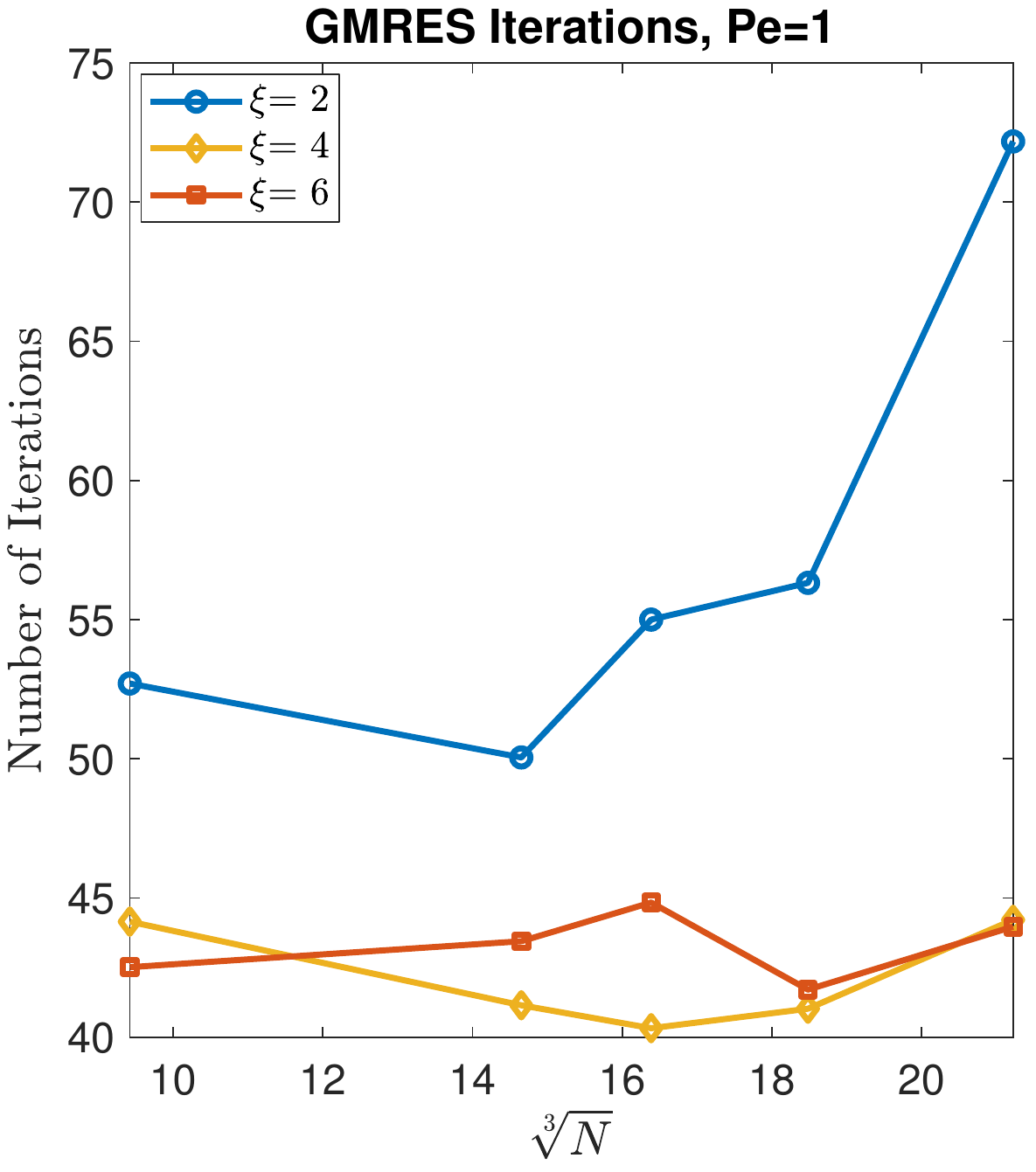}
	\label{fig:res3d_21}	
}
\subfloat[]
{
	\includegraphics[scale=0.5]{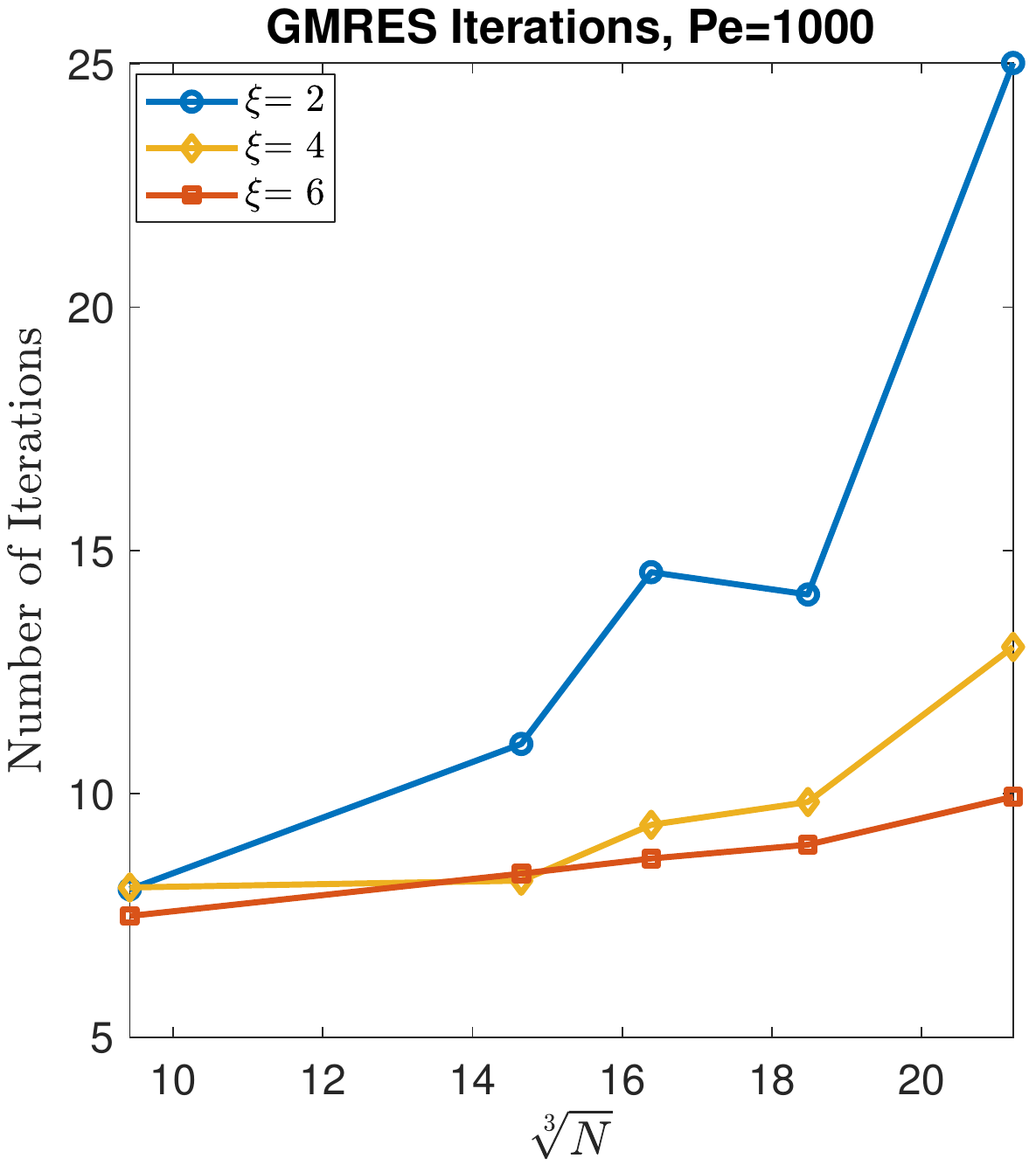}
	\label{fig:res3d_22}	
}
\caption{Top row: relative $\ell_2$ error vs $\sqrt[3]{N}$ as a function of approximation order $\xi$ for forced advection-diffusion in the unit ball. The dashed lines are lines of best fit indicating the slope (and hence convergence rate). Bottom row: average number of GMRES iterations per time-step as a function of $\sqrt[3]{N}$ and $\xi$.}
\label{fig:res3d}
\end{figure}
Next, we do a convergence study on the forced advection-diffusion equation in an irregular time-varying 3D domain. The initial domain is the \textbf{unit ball} with an embedded ball $E_3$ of radius $0.2$ centered at $(0.1, 0.2, 0.3)$. The initial simulation domain is $\Omega(t_0) = \{ \mathbb{B}^3 \setminus E_3\}$. We use 200 seed nodes on the boundary of $E_3$, and reconstruct the moving surface from these seed nodes using the (parametric) spherical, degree-8, polyharmonic spline RBF~\cite{SFKSISC2018}. In this case, our manufactured solution is
\begin{align}
c(\vx,t) = 1 + \sin(\pi x) \cos(\pi y) \cos(\pi z) \sin(\pi t),
\end{align}
and the incompressible velocity field $\vu(\vx,t) = (u,v,w)$ is given by
\begin{align}
\vu(\vx,t) = \sin\lf(\pi \|\vx\|_2^2\rt)\sin(\pi t)[yz,-2xz,xy].
\end{align}
The boundary conditions, time-steps, and Peclet numbers are chosen as in the 2D case, and the seed nodes are once again moved with the local fluid velocity, causing deformation while enforcing a no-slip condition. We simulate the PDE to time $t=0.5$, and measure errors against the manufactured solution. 

The results are shown in Figure \ref{fig:res3d}, plotted as function of $\sqrt[3]{N}$  (proportional to $1/h$). Once again, from Figures \ref{fig:res3d}(a) and \ref{fig:res3d}(b), we see that our results match the predicted spatial convergence rate of $h^{\xi}$ except for the slightly erratic convergence at $\rm{Pe} = 1000$ for $\xi=2$ (likely due to an insufficiently small timestep on the finer node set). Figures \ref{fig:res3d}(c) and \ref{fig:res3d}(d) show the average number of GMRES iterations per timestep for each of those Peclet numbers as a function of $\sqrt[3]{N}$. In both cases, increasing $\xi$ decreases the number of iterations (albeit erratically at low Peclet number). In addition, the number of iterations increase more slowly with $\sqrt[3]{N}$ as $\xi$ is increased, demonstrating the efficiency of higher order methods. As in the 2D case, Figure \ref{fig:res3d}(d) shows that a higher Peclet number requires fewer GMRES iterations to obtain the same tolerance, once again because of the improved conditioning due to smaller values of the diffusion coefficient $\nu$.

%% file: timings.tex
\subsection{Timings}
\label{sec:timings}
\begin{figure}[h!]
\centering
\subfloat[]
{
	\includegraphics[scale=0.5]{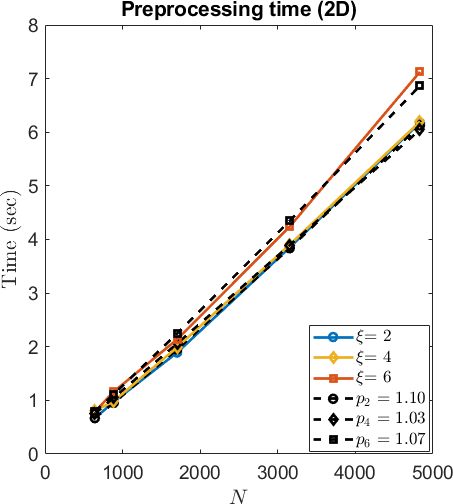}	
}
\subfloat[]
{
	\includegraphics[scale=0.5]{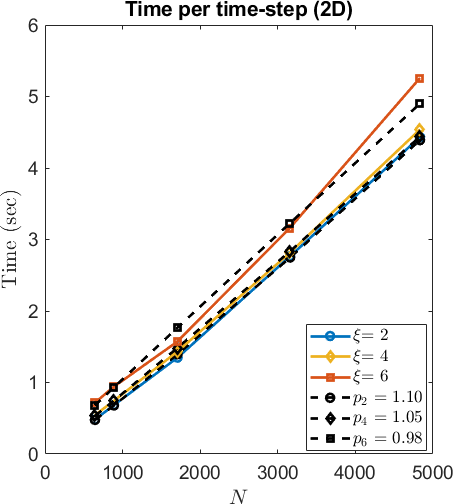}	
}

\subfloat[]
{
	\includegraphics[scale=0.5]{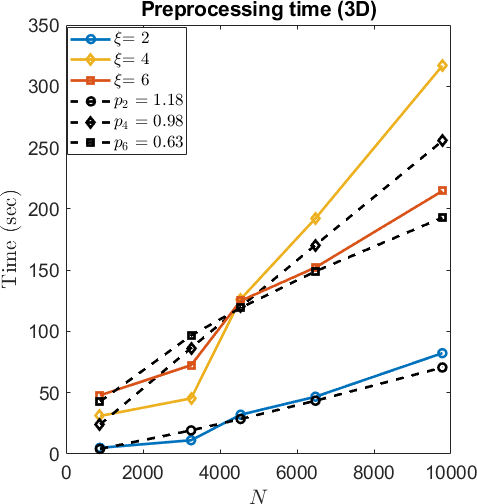}	
}
\subfloat[]
{
	\includegraphics[scale=0.5]{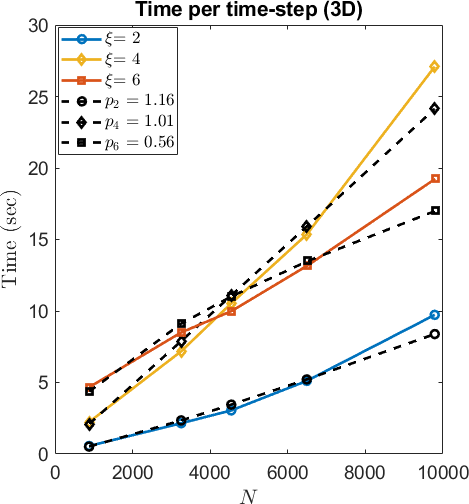}	
}
\caption{Top row: wall-clock time vs $N$ as a function of approximation order $\xi$ for forced advection-diffusion on the unit disk. The dashed lines are lines of best fit indicating the slope (and hence computational complexity). Bottom row: wall-clock time vs $N$ as a function of order $\xi$ for forced advection-diffusion in the unit ball. All timings are for Peclet number $\rm{Pe} = 1000$.}
\label{fig:time1}
\end{figure}
\begin{figure}[h!]
\centering
\subfloat[]
{
	\includegraphics[scale=0.5]{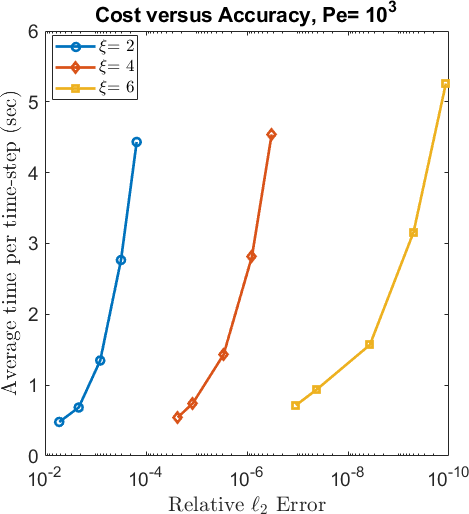}	
}
\subfloat[]
{
	\includegraphics[scale=0.5]{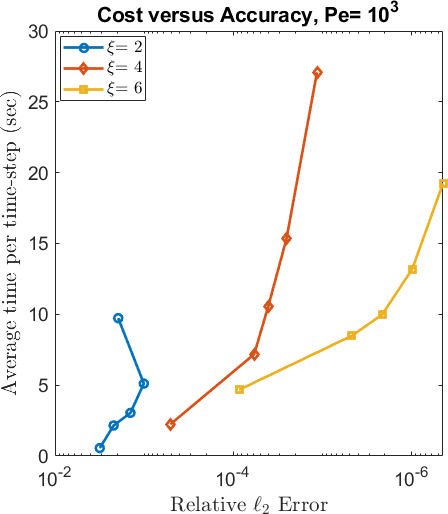}	
}
\caption{Cost vs accuracy as a function of approximation order $\xi$ for the 2D problem (left) and the 3D problem (right). The figure shows average time per time-step and the relative $\ell_2$ error as $N$ is increased for a given value of $\xi$. All timings are for Peclet number $\rm{Pe} = 1000$.}
\label{fig:time2}
\end{figure}
Next, we present timing results to verify the complexity estimates of our method. We present timings only for $\rm{Pe} = 1000$. This is because our goal is to verify the analysis from Section \ref{sec:comp_analysis}, which only abstractly includes the cost of the GMRES method. Since the results thus far have already shown that the number of GMRES iterations is smaller for $\rm{Pe} = 1000$, this case will be more ideal for verifying complexity. We present two types of results. All timings were run on a PC with an 8-core Intel i7-9700K CPU  clocked at a base speed of 3.6 GHz which had 32 GB of 2.6GHz RAM. Timings were accomplished in Matlab R2020b without any explicit parallelization of the code. 

In Figure \ref{fig:time1}, we present wall-clock times as a function of the number of nodes $N$ and the approximation order $\xi$ for both the 2D and 3D test problems. We separate the timings into preprocessing time (Figures \ref{fig:time1}(a) and (c)) and average per-step time (Figures \ref{fig:time1}(b) and (d)). The lines of best fit indicate that the slopes are close to 1, which indicates linear or quasi-linear complexity. In addition, it is clear that the cost per time-step is much lower than the preprocessing cost due to our strategy of copying as many RBF-FD weights as possible. The difference in costs is seen particularly in Figures \ref{fig:time1}(c) and (d) (for the 3D simulation), which show a full order of magnitude difference in time between a single preprocessing step and the average time-step. Even more interestingly, for the 3D results, we see that the sixth-order method ($\xi = 6$) is less expensive than the fourth-order method ($\xi = 4$) not just in terms of efficiency but in terms of wall-clock time. \revone{This is likely because the sixth-order method allows for a greater number of stable weights to be computed and retained per stencil, consequently allowing fewer stencils and recomputations in total.}

In Figure \ref{fig:time2}, we present the cost-accuracy tradeoffs for the 2D and 3D tests, again as a function of number of nodes $N$ and $\xi$. Here, we focus purely on the per-step costs, since the preprocessing costs are not important for long-running simulations. The figures all show that the higher order methods deliver greater accuracy for the same wall-clock time (and a given value of $N$), since our manufactured solutions are smooth. In addition, we once again see in Figure \ref{fig:time2}b that the sixth-order method takes even less time than the fourth-order method for the largest value of $N$ (top-most points of the curves) while delivering greater accuracy (in 3D). Both these timing tests confirm the quasilinear complexity of our method while also demonstrating the benefits of using high-order methods especially in higher dimensions. These tests also clearly demonstrate the importance of selective updates to RBF-FD weights and stencils using Algorithm \ref{alg:op_updates}.

%% file: coupled.tex
\section{A coupled problem}
\label{sec:application}
To fully demonstrate the utility of our numerical methods, we now apply them to solving a 3D coupled problem. In this problem, we track a chemical concentration $c(\vx,t)$ in a fluid inside the same \revone{initial} domain as in Section \ref{sec:res3d}, with the embedded boundary passively advected by an incompressible velocity field $\vu(\vx,t)$. The concentration $c(\vx,t)$ therefore satisfies the advection-diffusion equation (written using the Lagrangian/material derivative):
\begin{align}
\frac{d c}{dt} &= \nu \Delta c + f_1(\vx,t), \vx \in \Omega(t),
\label{eq:coupled_c}
\end{align}
where $f_1(\vx,t)$ is some forcing term, and $\Omega(t)$ is the time-varying irregular domain defined by the outer spherical boundary and the deforming inner embedded boundary $\Gamma(t)$. On the boundary of the unit ball, $c$ satisfies a time-varying inhomogeneous Neumann boundary condition:
\begin{align}
-\nu \frac{\partial c(\vx,t)}{\partial \vn} = g(\vx,t), \vx \in \mathbb{S},
\end{align}
which corresponds to an \emph{inward} flux of $c$. The inner moving embedded boundary $\Gamma(t)$ is viewed as a reactive ``zone'' on which chemicals can bind, unbind, and participate in other reactions. We track the \emph{bound} chemical surface density separately, and label it $C_B$. We assume that at each point on the embedded boundary $\Gamma(t)$, $C_B$ satisfies the following PDE:
\begin{align}
\frac{\partial C_B}{\partial t} + \vu \cdot \nabla C_B -  C_B \nabla_{\Gamma} \cdot \vu = k_{on}(C^{Tot} - C_B) - k_{off}C_B + k_{self}C_B(C^{Tot} - C_B) + f_2,
\label{eq:coupled_cb}
\end{align}
where $k_{on}$ and $k_{off}$ are the binding and unbinding rates of $C_B$, $k_{self}$ is the rate at which $C_B$ reacts with itself, $C^{Tot}$ is the total density of binding sites at each point of the reactive zone, $f_2$ is some forcing term, and $\nabla_{\Gamma} \cdot$ is the surface divergence operator. The last term on the left accounts for local surface changes due to the flow. As before, the first two terms in \eqref{eq:coupled_cb} can be combined into the material derivative, yielding the following equation for $C_B$:
\begin{align}
\frac{d C_B}{d t} -  C_B \nabla_{\Gamma} \cdot \vu = k_{on}(C^{Tot} - C_B) - k_{off}C_B + k_{self}C_B(C^{Tot} - C_B) + f_2,
\label{eq:coupled_cb2}
\end{align}
 Balancing fluxes at the interface $\Gamma(t)$ yields the following time-varying Robin boundary condition on $c(\vx,t)$:
\begin{align}
-\nu \frac{\partial c(\vx,t)}{\partial \vn} =- k_{on}(C^{Tot} - C_B) + k_{off}C_B, \vx \in \Gamma(t).
\label{eq:robin_bc}
\end{align}
This model represents a \emph{one-way} coupled bulk-surface problem, and is loosely inspired by similar problems arising in the context of platelet aggregation and coagulation. 

\subsection{A manufactured solution to the coupled problem}
We specify the concentration $c(\vx,t)$ to be
\begin{align}
c(\vx,t) = c(x,y,z,t) = 1 + \sin(\pi x) \cos(\pi y) \sin(\pi z) \sin(\pi t).
\end{align}
In addition, we \emph{prescribe} an incompressible velocity field $\vu(\vx,t) = (u,v,w)$ in the unit ball as
\begin{align}
\vu(\vx,t) = \sin(\pi t) \sin\lf(\pi \|\vx\|_2^2\rt) [yz,-2xz,xy].
\end{align}
Because points on the embedded boundary $\Gamma(t)$ will advect in this velocity field, the no-slip boundary is automatically satisfied on its surface. Then, we set \revone{the} forcing term $f_1(\vx,t)$ as
\begin{align}
f_1(\vx,t) = \frac{\partial c}{\partial t} + \vu \cdot \nabla c - \nu \Delta c.
\end{align}
The boundary condition function $g(\vx,t)$ is then obtained by applying the Neumann operator $-\nu \frac{\partial}{\partial \vn}$ to $c(\vx,t)$. We substitute $c(\vx,t)$ into \eqref{eq:robin_bc} and solve for $C_B$ to obtain
\begin{align}
C_B = \frac{-\nu \frac{\partial c}{\partial \vn} + k_{on}C^{Tot}}{k_{on} + k_{off}}.
\end{align}
Using $C_B$, we compute the forcing term $f_2$ as
\begin{align}
f_2 = \frac{d C_B}{d t} - C_B\nabla_{\Gamma} \cdot \vu - k_{on}\lf(C^{Tot}-C_B\rt)c_{amb} + k_{off}C_B - k_{self}C_B\lf(C^{Tot} - C_B\rt).
\end{align}
For a spatial convergence study, we define numerical errors by comparing against $c$ and $C_B$. Since $\vu$ is known analytically, the term $\nabla_{\Gamma}$ can be computed quasi-analytically using the numerically-computed normal vectors on $\Gamma(t)$. 

\subsection{Time-stepping for the coupled problem}
The fluid-phase chemicals $c(\vx,t)$ and the bound chemical \revone{density} $C_B$ are coupled through the boundary conditions in \eqref{eq:robin_bc}, which in turn reflect the binding and unbinding reactions in the ODE~\eqref{eq:coupled_cb}. To simulate this system efficiently, we use a simple time-splitting scheme that is a combination of SL updates for $c$ and full Lagrangian updates for $C_B$. Given $C_B^n$, $c^n$, and locations for the seed nodes on $\Gamma^n$, the algorithm is as follows:
\begin{enumerate}
\item Solve \eqref{eq:coupled_cb2} as an ODE using the semi-implicit BDF3 (SBDF3) method~\cite{Ascher97} to advance $C_B^n$ to $\tilde{C_B}$. The forcing term $f_2$ is treated implicitly in time, and all other terms are treated explicitly.
\item Advance the locations of the seed nodes on $\Gamma^n$ using $\vu(\vx,t)$ and the RK3 method. This \revone{advects} $\tilde{C_B}$ in a Lagrangian fashion to obtain $C_B^{n+1}$. Use the new seed node locations to construct a geometric representation of the embedded boundary and to interpolate $C_B$ to boundary nodes.
\item Use this representation and the interpolated $C_B^{n+1}$ to obtain the boundary conditions for $c^{n+1}$ on $\Gamma^{n+1}$ as:
\begin{align}
-\nu \frac{\partial c^{n+1}}{\partial \vn} = - k_{on}\lf(C^{Tot} - C_B^{n+1}\rt)+ k_{off} C_B^{n+1}. 
\end{align}
\item Update $c^n$ to $c^{n+1}$ using the SL method presented in Algorithm \ref{alg:sl_master}.
\end{enumerate}
It is important to note that a fully-coupled problem in which $c$ participates in the first term on the right hand side of \eqref{eq:coupled_cb2} would likely require high-order temporal splitting; we leave such an investigation to future work. For the purposes of this article, we have observed that the above scheme produces highly accurate results both for $c$ and $C_B$. It is also worth noting that if $\vu$ is not known analytically, the $\nabla_{\Gamma} \cdot \vu$ term must be computed numerically. This can also be done using a stabilized version of overlapped RBF-FD specialized to manifolds~\cite{SNKJCP2018,SNWSISC2020}.

\subsection{Results}
\begin{figure}[h!]
\centering
\subfloat[]
{
	\includegraphics[scale=0.5]{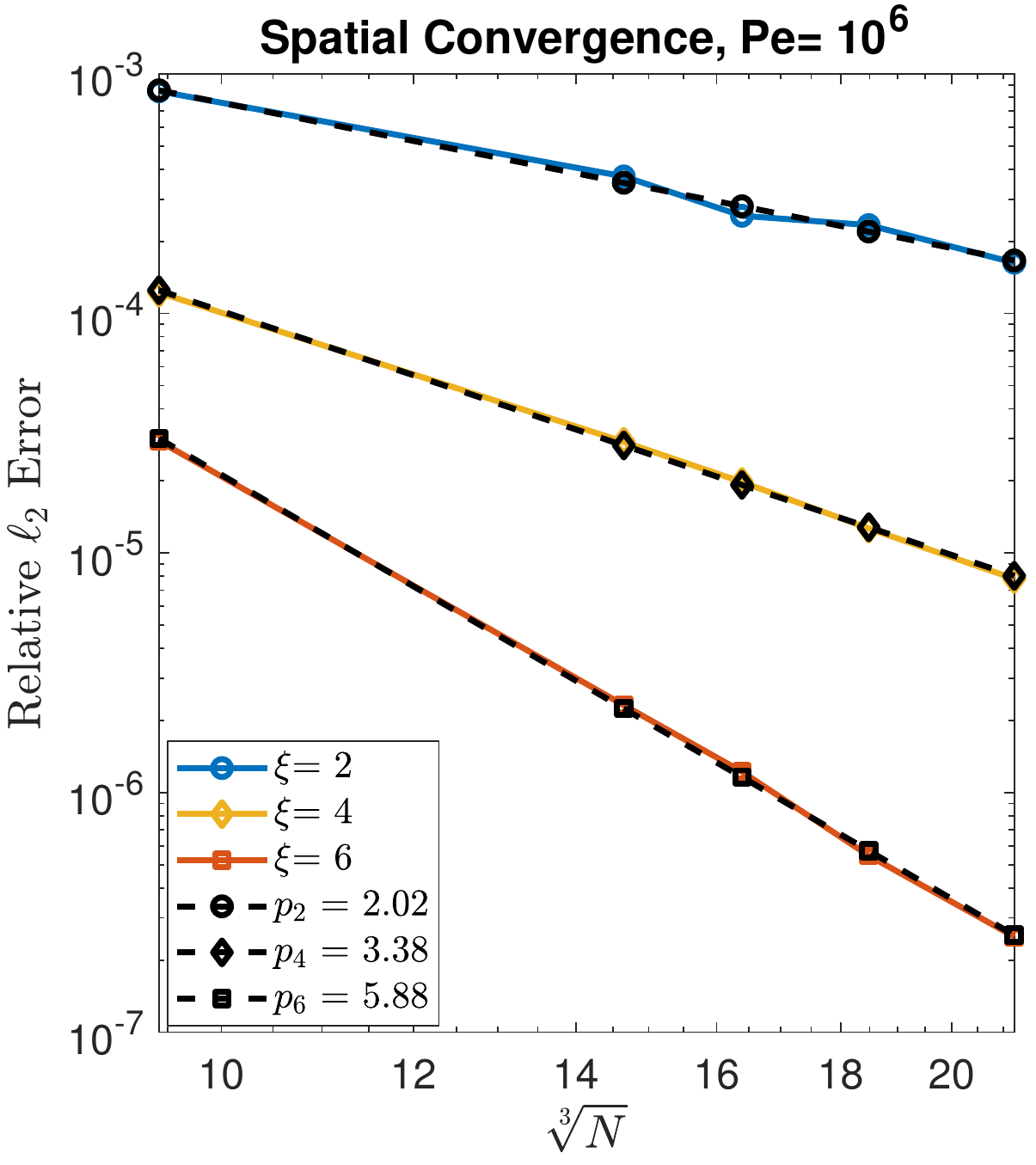}
	\label{fig:bs_11}	
}
\subfloat[]
{
	\includegraphics[scale=0.5]{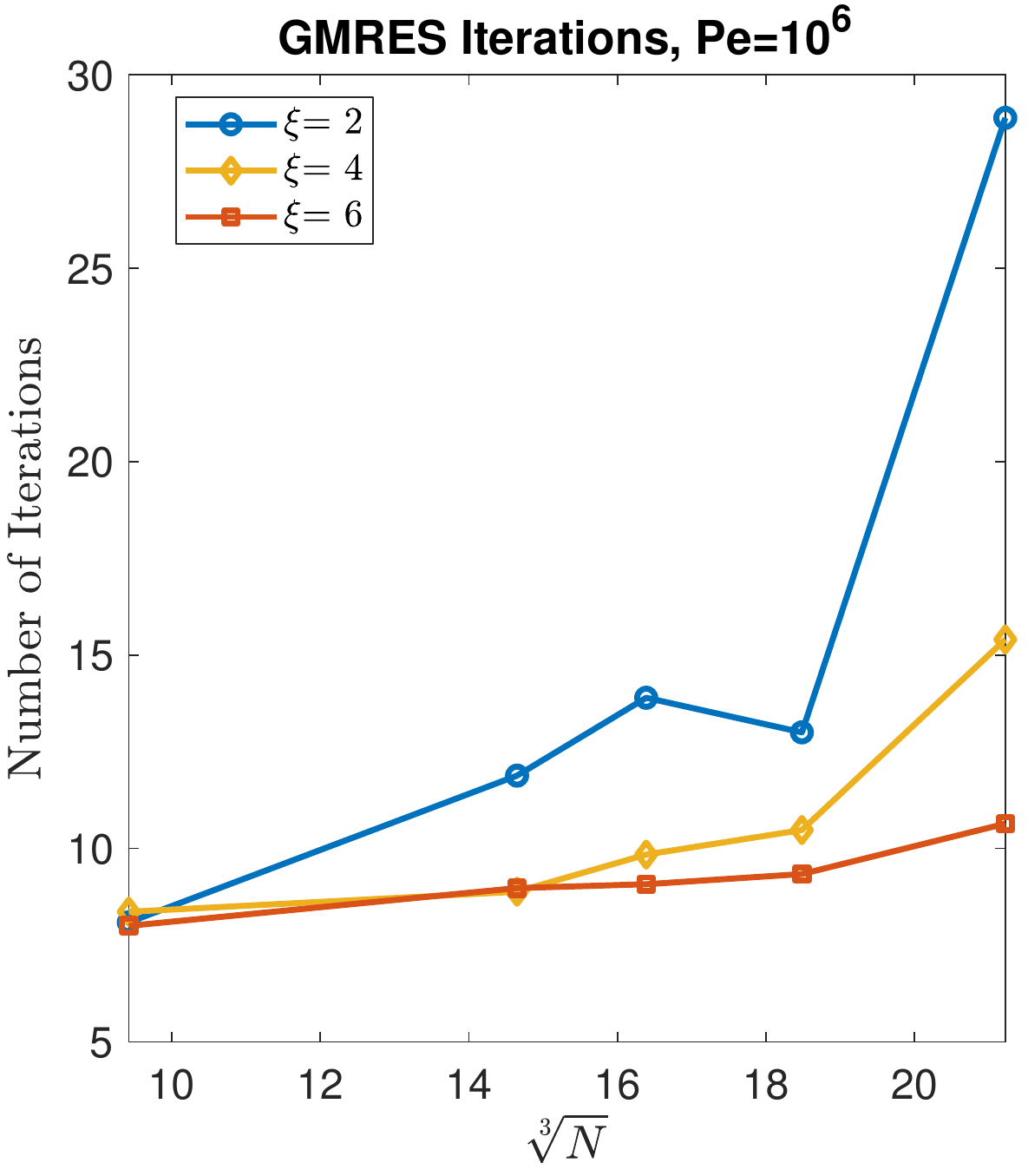}
	\label{fig:bs_21}	
}
\caption{Left: relative errors vs node spacing as a function of approximation order $\xi$. Right: number of GMRES iterations as a function of $\xi$.}
\label{fig:bs_conv_plots}
\end{figure}
To more fully test the stability of our method, we now set the diffusion coefficient to $\nu = 10^{-6}$, and scale the velocity field $\vu$ to obtain a Peclet number of $\rm{Pe} = 10^6$. We select the time-step as in the advection-diffusion test cases, run convergence studies to time $t=0.5$ using $\xi = 2,4$ and $6$, and measure both the relative error in $c$ and the average number of GMRES iterations. The results are shown in Figure \ref{fig:bs_conv_plots}. Much like our other tests, Figure \ref{fig:bs_conv_plots}(a) shows that the spatial error in the numerical approximation to $c$ decreases at the rate of approximately $h^{\xi}$. In addition, Figure \ref{fig:bs_conv_plots}(b) shows that the number of GMRES iterations increases approximately linearly with the node spacing, and decreases as $\xi$ is increased. It is also worth noting that the relative error in $C_B$ is purely temporal and is on the order of $10^{-8}$ for all the values of $\triangle t$ used in this test. This application clearly demonstrates the utility of our overall framework in solving coupled bulk-surface problems on moving domains.

%% file: Discussion.tex
\section{Summary and Future Work}
\label{sec:summary}

In this article, we presented a high-order numerical method for simulating the advection-diffusion equation on domains with time-varying embedded boundaries. Our method uses semi-Lagrangian (SL) advection combined with an Eulerian formulation for diffusion, both applied in the context of a rapid node adaptation algorithm. These techniques rely on a generalization of the overlapped RBF-FD method that replaces the hand-tuned overlap parameter $\delta$ with a pair of automatically-computed stability indicators. We also presented a novel automatic algorithm for updating RBF-FD interpolation stencils and differentiation matrices on a time-varying node set. We conducted an informal error analysis to show the high-order convergence rate of our method, and conducted a computational complexity analysis to determine the efficiency of our algorithms. We verified the aforementioned high-order convergence rates on both 2D and 3D advection-diffusion problems on irregular domains with moving embedded boundaries. In addition, we demonstrated high-order convergence rates on a more complicated 3D coupled problem.

While our informal error analysis provided estimates that were verified numerically, a formal error analysis of semi-Lagrangian RBF-FD methods for advection-diffusion equations is absent in the literature. We plan to address this issue in future work. In addition, in order to apply our methods to problems requiring very large node sets, the overlapped RBF-FD method (and the overall method presented in this article) must be parallelized for distributed and shared memory architectures. We plan to tackle this in future work as well.

%% file: article.bbl
\begin{thebibliography}{}

\bibitem[\protect\astroncite{Ascher et~al.}{1997}]{Ascher97}
Ascher, U.~M., Ruuth, S.~J., and Wetton, B. T.~R. (1997).
\newblock {Implicit-Explicit Methods For Time-Dependent PDEs}.
\newblock {\em SIAM J. Numer. Anal}, 32:797--823.

\bibitem[\protect\astroncite{Barnett}{2015}]{BarnettPHS}
Barnett, G.~A. (2015).
\newblock {\em A Robust RBF-FD Formulation based on Polyharmonic Splines and
  Polynomials}.
\newblock PhD thesis, University of Colorado Boulder.

\bibitem[\protect\astroncite{Bayona}{2019}]{Bayona2019}
Bayona, V. (2019).
\newblock Comparison of moving least squares and {RBF+ poly} for interpolation
  and derivative approximation.
\newblock {\em Journal of Scientific Computing}, 81(1):486--512.

\bibitem[\protect\astroncite{Bayona et~al.}{2019}]{BayonaBoundary}
Bayona, V., Flyer, N., and Fornberg, B. (2019).
\newblock On the role of polynomials in rbf-fd approximations: Iii. behavior
  near domain boundaries.
\newblock {\em Journal of Computational Physics}, 380:378--399.

\bibitem[\protect\astroncite{Bayona et~al.}{2017}]{FlyerElliptic}
Bayona, V., Flyer, N., Fornberg, B., and Barnett, G.~A. (2017).
\newblock On the role of polynomials in {RBF-FD} approximations: {II}.
  {N}umerical solution of elliptic {PDE}s.
\newblock {\em J. Comput. Phys.}, 332:257--273.

\bibitem[\protect\astroncite{Bayona et~al.}{2010}]{Bayona2010}
Bayona, V., Moscoso, M., Carretero, M., and Kindelan, M. (2010).
\newblock {RBF-FD} formulas and convergence properties.
\newblock {\em J. Comput. Phys.}, 229(22):8281--8295.

\bibitem[\protect\astroncite{Behrens and Iske}{2002}]{iske2002}
Behrens, J. and Iske, A. (2002).
\newblock {Grid-free adaptive semi-{L}agrangian advection using radial basis
  functions}.
\newblock {\em Comput. Math. Appl.}, 43(3):319--327.

\bibitem[\protect\astroncite{Benzi et~al.}{2005}]{Benzi2005}
Benzi, M., Golub, G.~H., and Liesen, J. (2005).
\newblock Numerical solution of saddle point problems.
\newblock {\em Acta Numerica}, 14:1--137.

\bibitem[\protect\astroncite{Benzi and Wathen}{2008}]{Benzi2008}
Benzi, M. and Wathen, A.~J. (2008).
\newblock Some preconditioning techniques for saddle point problems.
\newblock In {\em Model order reduction: theory, research aspects and
  applications}, pages 195--211. Springer.

\bibitem[\protect\astroncite{Bonaventura and Ferretti}{2014}]{BVSDE1}
Bonaventura, L. and Ferretti, R. (2014).
\newblock Semi-{L}agrangian methods for parabolic problems in divergence form.
\newblock {\em SIAM Journal on Scientific Computing}, 36(5):A2458--A2477.

\bibitem[\protect\astroncite{Bonaventura and Ferretti}{2016}]{BVSDE2}
Bonaventura, L. and Ferretti, R. (2016).
\newblock Flux form semi-{L}agrangian methods for parabolic problems.
\newblock {\em Communications in Applied and Industrial Mathematics},
  7(3):56--73.

\bibitem[\protect\astroncite{Bowman et~al.}{2015}]{LagRe}
Bowman, J.~C., Yassaei, M.~A., and Basu, A. (2015).
\newblock A fully {L}agrangian advection scheme.
\newblock {\em Journal of Scientific Computing}, 64(1):151--177.

\bibitem[\protect\astroncite{Calhoun and LeVeque}{2000}]{Calhoun99acartesian}
Calhoun, D. and LeVeque, R.~J. (2000).
\newblock A {C}artesian grid finite-volume method for the advection-diffusion
  equation in irregular geometries.
\newblock {\em J. Comput. Phys.}, 157(1):143--180.

\bibitem[\protect\astroncite{Colonius and Taira}{2008}]{ColoniusTaira08}
Colonius, T. and Taira, K. (2008).
\newblock A fast immersed boundary method using a nullspace approach and
  multi-domain far-field boundary conditions.
\newblock {\em Compu. Methods Appl. Mech. Engrg.}, 197:2131--2146.

\bibitem[\protect\astroncite{Davydov and Oanh}{2011}]{Davydov2011}
Davydov, O. and Oanh, D.~T. (2011).
\newblock Adaptive meshless centres and {RBF} stencils for {P}oisson equation.
\newblock {\em J. Comput. Phys.}, 230(2):287--304.

\bibitem[\protect\astroncite{Davydov and Schaback}{2018}]{DavydovSchaback2018}
Davydov, O. and Schaback, R. (2018).
\newblock Minimal numerical differentiation formulas.
\newblock {\em Numerische Mathematik}, 140(3):555--592.

\bibitem[\protect\astroncite{Degond and Mas-Gallic}{1989}]{PSE1}
Degond, P. and Mas-Gallic, S. (1989).
\newblock The weighted particle method for convection-diffusion equations. part
  1: The case of an isotropic viscosity.
\newblock {\em Mathematics of Computation}, 53(188):485--507.

\bibitem[\protect\astroncite{Duff and Koster}{2001}]{MatlabEquilibrate}
Duff, I.~S. and Koster, J. (2001).
\newblock On algorithms for permuting large entries to the diagonal of a sparse
  matrix.
\newblock {\em SIAM Journal on Matrix Analysis and Applications},
  22(4):973--996.

\bibitem[\protect\astroncite{Fadlun et~al.}{2000}]{IB-forcing-Fadlun-JCP2000}
Fadlun, E., Verzicco, R., Orlandi, P., and Mohd-Yusof, J. (2000).
\newblock Combined immersed-boundary finite-difference methods for
  three-dimensional complex flow simulations.
\newblock {\em J. Comput. Phys.}, 161:35--60.

\bibitem[\protect\astroncite{Falcone and Ferretti}{1998}]{FalconeFerretti}
Falcone, M. and Ferretti, R. (1998).
\newblock Convergence analysis for a class of high-order semi-{L}agrangian
  advection schemes.
\newblock {\em SIAM Journal on Numerical Analysis}, 35(3):909--940.

\bibitem[\protect\astroncite{Falcone and Ferretti}{2013}]{FFSL}
Falcone, M. and Ferretti, R. (2013).
\newblock {\em Semi-{L}agrangian approximation schemes for linear and
  Hamilton—Jacobi equations}.
\newblock SIAM.

\bibitem[\protect\astroncite{Fasshauer}{2007}]{Fasshauer:2007}
Fasshauer, G.~E. (2007).
\newblock {\em {Meshfree Approximation Methods with {MATLAB}}}.
\newblock Interdisciplinary Mathematical Sciences - Vol. 6. World Scientific
  Publishers, Singapore.

\bibitem[\protect\astroncite{Flyer et~al.}{2016a}]{FlyerNS}
Flyer, N., Barnett, G.~A., and Wicker, L.~J. (2016a).
\newblock Enhancing finite differences with radial basis functions: Experiments
  on the {N}avier-{S}tokes equations.
\newblock {\em J. Comput. Phys.}, 316:39--62.

\bibitem[\protect\astroncite{Flyer et~al.}{2016b}]{FlyerPHS}
Flyer, N., Fornberg, B., Bayona, V., and Barnett, G.~A. (2016b).
\newblock {On the role of polynomials in RBF-FD approximations: I.
  Interpolation and accuracy}.
\newblock {\em J. Comput. Phys.}, 321:21--38.

\bibitem[\protect\astroncite{Flyer et~al.}{2012}]{FlyerLehto2012}
Flyer, N., Lehto, E., Blaise, S., Wright, G.~B., and St-Cyr, A. (2012).
\newblock {A guide to {RBF}-generated finite differences for nonlinear
  transport: shallow water simulations on a sphere}.
\newblock {\em J. Comput. Phys.}, 231:4078--4095.

\bibitem[\protect\astroncite{Flyer and Wright}{2007}]{FlyerWright:2007}
Flyer, N. and Wright, G.~B. (2007).
\newblock Transport schemes on a sphere using radial basis functions.
\newblock {\em J. Comput. Phys.}, 226:1059--1084.

\bibitem[\protect\astroncite{Flyer and Wright}{2009}]{FlyerWright:2009}
Flyer, N. and Wright, G.~B. (2009).
\newblock {A radial basis function method for the shallow water equations on a
  sphere}.
\newblock {\em Proc. Roy. Soc. A}, 465:1949--1976.

\bibitem[\protect\astroncite{Fornberg and Lehto}{2011}]{FoL11}
Fornberg, B. and Lehto, E. (2011).
\newblock {Stabilization of {RBF}-generated finite difference methods for
  convective {PDE}s}.
\newblock {\em J. Comput. Phys.}, 230:2270--2285.

\bibitem[\protect\astroncite{Fuselier and Wright}{2013}]{FuselierWright2013}
Fuselier, E.~J. and Wright, G.~B. (2013).
\newblock A high-order kernel method for diffusion and reaction-diffusion
  equations on surfaces.
\newblock {\em J. Sci. Comput.}, 56(3):535--565.

\bibitem[\protect\astroncite{Glowinski et~al.}{1998}]{Glowinski}
Glowinski, R., Pan, T., and P\'{e}riaux, J. (1998).
\newblock Distributed {L}agrange multiplier methods for incompressible viscous
  flow around moving rigid bodies.
\newblock {\em Comput. Methods Appl. Mech. Engrg.}, 151:181--194.

\bibitem[\protect\astroncite{Goldstein
  et~al.}{1993}]{IB-externalforce-Goldstein-JCP1993}
Goldstein, D., Handler, R., and Sirovich, L. (1993).
\newblock Modeling a no-slip flow boundary with an external force field.
\newblock {\em J. Comput. Phys.}, 105:354--366.

\bibitem[\protect\astroncite{Gritton et~al.}{2017}]{Gritton2017}
Gritton, C., Guilkey, J., Hooper, J., Bedrov, D., Kirby, R.~M., and Berzins, M.
  (2017).
\newblock Using the material point method to model chemical/mechanical coupling
  in the deformation of a silicon anode.
\newblock {\em Modelling and Simulation in Materials Science and Engineering},
  25(4):045005.

\bibitem[\protect\astroncite{Johansen and Colella}{1998}]{EB-poisson-JCP1998}
Johansen, H.~S. and Colella, P. (1998).
\newblock A {C}artesian grid embedded boundary method for poisson's equation on
  irregular domains.
\newblock {\em J. Comput. Phys.}, 147:60--85.

\bibitem[\protect\astroncite{Kim
  et~al.}{2001}]{IB-forcing_sink-KimChoi-JCP2001}
Kim, J., Kim, D., and Choi, H. (2001).
\newblock An immersed-boundary finite-volume method for simulations of flow in
  complex geometries.
\newblock {\em J. Comput. Phys.}, 171:132--150.

\bibitem[\protect\astroncite{Le~Roux et~al.}{1997}]{TELA:TELA0009}
Le~Roux, D.~Y., Lin, C.~A., and Staniforth, A. (1997).
\newblock An accurate interpolating scheme for semi-{L}agrangian advection on
  an unstructured mesh for ocean modelling.
\newblock {\em Tellus A}, 49(1):119--138.

\bibitem[\protect\astroncite{Lehto et~al.}{2017}]{LSWSISC2017}
Lehto, E., Shankar, V., and Wright, G.~B. (2017).
\newblock A radial basis function ({RBF}) compact finite difference ({FD})
  scheme for reaction-diffusion equations on surfaces.
\newblock {\em SIAM J. Sci. Comput.}, 39:A2129--A2151.

\bibitem[\protect\astroncite{Leiderman and Fogelson}{2014}]{LEIDERMAN:2014:OMM}
Leiderman, K. and Fogelson, A.~L. (2014).
\newblock {An Overview of Mathematical Modeling of Thrombus Formation Under
  Flow}.
\newblock {\em Thromb. Res.}, 133 Suppl:S12--S14.

\bibitem[\protect\astroncite{Leiderman and Fogelson}{2011}]{LEIDERMAN:2011:GWF}
Leiderman, K.~M. and Fogelson, A.~L. (2011).
\newblock {Grow with the flow: a spatial-temporal model of platelet deposition
  and blood coagulation under flow}.
\newblock {\em Math. Med. Biol.}, 28:47--84.

\bibitem[\protect\astroncite{Leiderman and Fogelson}{2013}]{LEIDERMAN:2013:IHT}
Leiderman, K.~M. and Fogelson, A.~L. (2013).
\newblock The influence of hindered transport on the development of platelet
  thrombi under flow.
\newblock {\em Bull. of Math. Biol.}, 75:1255--1283.

\bibitem[\protect\astroncite{LeVeque and Li}{1994}]{IIM1}
LeVeque, R.~J. and Li, Z. (1994).
\newblock The immersed interface method for elliptic equations with
  discontinuous coefficients and singular sources.
\newblock {\em SIAM J. Numer. Anal.}, 31:1001--1025.

\bibitem[\protect\astroncite{McCorquodale
  et~al.}{2001}]{McCorquodaleColellaJohansen01}
McCorquodale, P., Colella, P., and Johansen, H. (2001).
\newblock A {C}artesian grid embedded boundary method for the heat equation on
  irregular domains.
\newblock {\em J. Comput. Phys.}, 173:620--635.

\bibitem[\protect\astroncite{Mohd-Yusof}{1997}]{ForcingBC-Mohd-Yusof1997}
Mohd-Yusof, J. (1997).
\newblock Combined immersed-boundary/{B}-spline methods for simulations of flow
  in complex geometries.
\newblock {\em Annu. Res.Briefs, Cent. Turbul. Res.}, pages 317--328.

\bibitem[\protect\astroncite{Peskin}{1972}]{IBM1}
Peskin, C.~S. (1972).
\newblock Flow pattern around heart valves: a numerical method.
\newblock {\em J. Comput. Phys.}, 10:252--271.

\bibitem[\protect\astroncite{Peskin}{1977}]{IBM2}
Peskin, C.~S. (1977).
\newblock Numerical analysis of blood flow in the heart.
\newblock {\em J. Comput. Phys}, 25:220--252.

\bibitem[\protect\astroncite{Peskin}{2002}]{PESKIN:2002:IBM}
Peskin, C.~S. (2002).
\newblock The immersed boundary method.
\newblock {\em Acta Numerica}, 11:479--517.

\bibitem[\protect\astroncite{Piret}{2012}]{Piret2012}
Piret, C. (2012).
\newblock {The orthogonal gradients method: A radial basis functions method for
  solving partial differential equations on arbitrary surfaces}.
\newblock {\em J. Comput. Phys.}, 231(20):4662--4675.

\bibitem[\protect\astroncite{Piret and Dunn}{2016}]{Piret2016}
Piret, C. and Dunn, J. (2016).
\newblock Fast {RBF} {OGr} for solving pdes on arbitrary surfaces.
\newblock {\em AIP Conference Proceedings}, 1776(1).

\bibitem[\protect\astroncite{Saad}{2003}]{Saad2003}
Saad, Y. (2003).
\newblock {\em Iterative methods for sparse linear systems}.
\newblock SIAM.

\bibitem[\protect\astroncite{Shankar}{2017}]{ShankarJCP2017}
Shankar, V. (2017).
\newblock The overlapped radial basis function-finite difference ({RBF-FD})
  method: A generalization of {RBF-FD}.
\newblock {\em J. Comput. Phys.}, 342:211--228.

\bibitem[\protect\astroncite{Shankar and Fogelson}{2018}]{SFJCP2018}
Shankar, V. and Fogelson, A.~L. (2018).
\newblock Hyperviscosity-based stabilization for radial basis function-finite
  difference (rbf-fd) discretizations of {advection– diffusion} equations.
\newblock {\em J. Comput. Phys.}, 372:616 -- 639.

\bibitem[\protect\astroncite{Shankar et~al.}{2018a}]{SFKSISC2018}
Shankar, V., Kirby, R., and Fogelson, A. (2018a).
\newblock Robust node generation for mesh-free discretizations on irregular
  domains and surfaces.
\newblock {\em SIAM Journal on Scientific Computing}, 40(4):A2584--A2608.

\bibitem[\protect\astroncite{Shankar et~al.}{2018b}]{SNKJCP2018}
Shankar, V., Narayan, A., and Kirby, R.~M. (2018b).
\newblock Rbf-loi: Augmenting radial basis functions (rbfs) with least
  orthogonal interpolation (loi) for solving pdes on surfaces.
\newblock {\em Journal of Computational Physics}, 373:722--735.

\bibitem[\protect\astroncite{Shankar and Wright}{2018}]{SWJCP2018}
Shankar, V. and Wright, G.~B. (2018).
\newblock Mesh-free semi-lagrangian methods for transport on a sphere using
  radial basis functions.
\newblock {\em J. Comput. Phys.}, 366(C):170--190.

\bibitem[\protect\astroncite{Shankar et~al.}{2014a}]{SWFKIJNMF2014}
Shankar, V., Wright, G.~B., Fogelson, A.~L., and Kirby, R.~M. (2014a).
\newblock A radial basis function ({RBF}) finite difference method for the
  simulation of reaction-diffusion equations on stationary platelets within the
  augmented forcing method.
\newblock {\em Inter. J. Numer. Methods Fluids}, 75(1):1--22.

\bibitem[\protect\astroncite{Shankar et~al.}{2014b}]{SWFKJSC2014}
Shankar, V., Wright, G.~B., Kirby, R.~M., and Fogelson, A.~L. (2014b).
\newblock A radial basis function ({RBF})-finite difference ({FD}) method for
  diffusion and reaction--diffusion equations on surfaces.
\newblock {\em J. Sci. Comput.}, 63(3):745--768.

\bibitem[\protect\astroncite{Shankar et~al.}{2015}]{SWFKIJNMF2015}
Shankar, V., Wright, G.~B., Kirby, R.~M., and Fogelson, A.~L. (2015).
\newblock Augmenting the immersed boundary method with radial basis functions
  (rbfs) for the modeling of platelets in hemodynamic flows.
\newblock {\em International Journal for Numerical Methods in Fluids},
  79(10):536--557.

\bibitem[\protect\astroncite{Shankar et~al.}{2020}]{SNWSISC2020}
Shankar, V., Wright, G.~B., and Narayan, A. (2020).
\newblock A robust hyperviscosity formulation for stable {RBF-FD}
  discretizations of {Advection-Diffusion-Reaction} equations on manifolds.
\newblock {\em SIAM Journal on Scientific Computing}, 42(4):A2371--A2401.

\bibitem[\protect\astroncite{Smolarkiewicz and
  Margolin}{1997}]{SmolarkiewiczMargolin97}
Smolarkiewicz, P.~K. and Margolin, L.~G. (1997).
\newblock On forward-in-time differencing for fluids: an
  {E}ulerian/semi-{L}agrangian non-hydrostatic model for stratified flows.
\newblock {\em Atmos.-Ocean}, 35:127--152.

\bibitem[\protect\astroncite{Smolarkiewicz and Pudykiewicz}{1992}]{SmolPudy92}
Smolarkiewicz, P.~K. and Pudykiewicz, J.~A. (1992).
\newblock A class of semi-{L}agrangian approximations for fluids.
\newblock {\em J. Atmos. Sci.}, 49(22):2082--2096.

\bibitem[\protect\astroncite{Staniforth and
  C{\^{o}}t{\'{e}}}{1991}]{staniforth1991SL}
Staniforth, A. and C{\^{o}}t{\'{e}}, J. (1991).
\newblock {Semi-{L}agrangian integration schemes for atmospheric models---a
  review}.
\newblock {\em Mon. Wea. Rev.}, 119(9):2206--2223.

\bibitem[\protect\astroncite{Staniforth and Wood}{2008}]{StaniforthWoodJCP2008}
Staniforth, A. and Wood, N. (2008).
\newblock Aspects of the dynamical core of a nonhydrostatic, deep-atmosphere,
  unified weather and climate-prediction model.
\newblock {\em J. Comput. Phys.}, 227(7):3445--3464.

\bibitem[\protect\astroncite{Stein et~al.}{2016}]{IBSE1}
Stein, D.~B., Guy, R.~D., and Thomases, B. (2016).
\newblock Immersed boundary smooth extension: a high-order method for solving
  {PDE} on arbitrary smooth domains using fourier spectral methods.
\newblock {\em J. Comput. Phys.}, 304:252--274.

\bibitem[\protect\astroncite{Stein et~al.}{2017}]{IBSE2}
Stein, D.~B., Guy, R.~D., and Thomases, B. (2017).
\newblock Immersed boundary smooth extension ({IBSE}): A high-order method for
  solving incompressible flows in arbitrary smooth domains.
\newblock {\em J. Comput. Phys.}, 335:155--178.

\bibitem[\protect\astroncite{Trask et~al.}{2015}]{trask2015Scalable}
Trask, N., Maxey, M., Kim, K., Perego, M., Parks, M.~L., Yang, K., and Xu, J.
  (2015).
\newblock A scalable consistent second-order {SPH} solver for unsteady low
  {R}eynolds number flows.
\newblock {\em Computer Methods in Applied Mechanics and Engineering},
  289:155--178.

\bibitem[\protect\astroncite{Udaykumar et~al.}{1999}]{Udaykumar1999}
Udaykumar, H.~S., Mittal, R., and Shyy, W. (1999).
\newblock Computation of solid-liquid phase fronts in the sharp interface limit
  on fixed grids.
\newblock {\em J. Comput. Phys.}, 153:535--574.

\bibitem[\protect\astroncite{Wendland}{2005}]{Wendland:2004}
Wendland, H. (2005).
\newblock {\em {Scattered data approximation}}, volume~17 of {\em Cambridge
  Monogr. Appl. Comput. Math.}
\newblock Cambridge University Press, Cambridge.

\bibitem[\protect\astroncite{Wright and Fornberg}{2006}]{Wright200699}
Wright, G.~B. and Fornberg, B. (2006).
\newblock {Scattered node compact finite difference-type formulas generated
  from radial basis functions}.
\newblock {\em J. Comput. Phys.}, 212(1):99--123.

\bibitem[\protect\astroncite{Xiu and Karniadakis}{2001}]{Xiu2002}
Xiu, D. and Karniadakis, G.~E. (2001).
\newblock {A semi-{L}agrangian high-order method for {N}avier--{S}tokes
  equations}.
\newblock {\em J. Comput. Phys.}, 172(2):658--684.

\bibitem[\protect\astroncite{Yao and Fogelson}{2012}]{YaoFogelson2012}
Yao, L. and Fogelson, A.~L. (2012).
\newblock Simulations of chemical transport and reaction in a suspension of
  cells {I}: an augmented forcing point method for the stationary case.
\newblock {\em Inter. J. Numer. Methods Fluids}, 69(11):1736--1752.

\bibitem[\protect\astroncite{Ye et~al.}{1999}]{IB-cutcell-Ye-JCP1999}
Ye, T., Mittal, R., Udaykumar, H.~S., and Shyy, W. (1999).
\newblock An accurate {C}artesian grid method for viscous incompressible flows
  with complex immersed boundaries.
\newblock {\em J. Comput. Phys.}, 156:209--240.

\end{thebibliography}
